\documentclass[a4paper,12pt]{article}
\usepackage[utf8]{inputenc}
\usepackage{amssymb}
\usepackage{amsmath}
\usepackage{MnSymbol} 
\usepackage{cancel}
\usepackage{ulem}
\usepackage{mathtools}
\usepackage{blkarray} 
\usepackage{tikz}
\usetikzlibrary{arrows}
\usepackage{tikz-cd}
\usepackage{fancyhdr}
\usepackage{authblk}

\usepackage{xcolor}

\usepackage{hyperref}

\usepackage{caption}
\usepackage{subcaption}

\newtheorem{definition}{Definition}
\newtheorem{example}{Example}

\DeclareMathOperator{\Ima}{Im} 
\DeclareMathOperator{\Hom}{Hom} 
\newcommand{\dflcat}{\mathbf{dFlag}^{\hookrightarrow}}

\author[1]{Benjamin Jones \thanks{jones657@msu.edu}}
\author[1,2,3]{Guo-Wei Wei \thanks{weig@msu.edu}}
\affil[1]{Department of Mathematics, Michigan State University, MI, 48824, USA}
\affil[2]{Department of Electrical and Computer Engineering, Michigan State University, MI 48824, USA}
\affil[3]{Department of Biochemistry and Molecular Biology, Michigan State University, MI 48824, USA}

\makeatletter
    \renewcommand*{\@fnsymbol}[1]{\ensuremath{\ifcase#1\or \dagger\or *\or *\or
   \mathsection\or \else\@ctrerr\fi}}
\makeatother
\date{}

\title{Persistent Directed Flag Laplacian}
\date{December 2023}

\begin{document}
\maketitle
\paragraph{Abstract} 
Topological data analysis (TDA) has had enormous success in science and engineering in the past decade.  Persistent topological Laplacians (PTLs) overcome some limitations of persistent homology, a key technique in TDA, and  provide substantial insight to the behavior of various geometric and topological objects. 
This work extends PTLs to directed flag complexes, which are an exciting generalization to flag complexes, also known as clique complexes, that arise naturally in many situations. We introduce the directed flag Laplacian and show that the proposed persistent directed flag Laplacian (PDFL) is a distinct way of analyzing these flag complexes. 
Example calculations are provided to demonstrate the potential of the proposed PDFL in real world applications.  


\newpage
    \tableofcontents
\newpage

\section{Introduction}
Topological data analysis (TDA) is a rapidly growing field that uses tools from a multitude of disciplines to study connections between geometric, topological, directional, time-based, and functional information. TDA builds upon the fundamental objects of homology and cohomology through a filtration process, which gives rise to persistence modules.
Persistence modules, including persistent (co)homology, are some of the most studied objects in TDA, and can convey information about a system in different dimensions, scales of size, and/or scales of time \cite{zomorodian2004computing,carlsson2009topology}. They can be analyzed using persistent barcodes \cite{ghrist2008barcodes}, persistent diagrams, persistent images \cite{adams2017persistence},  persistence landscapes \cite{persistenceLandscapes2015Bubenik}, and persistent Betti numbers. These methods have been used to study protein folding \cite{xia2014persistent}, protein-ligand binding \cite{cang2017topologynet}, virus mutation \cite{chen2020mutations}, chemistry \cite{townsend2020representation},  dynamical systems \cite{pereaTopologicalTimeSeries2019}, signal processing \cite{myersTeaspoonComprehensivePython2020}. In fact,  persistent homology owns its success to  advanced  machine learning models \cite{henselSurveyTopologicalMachine2021}, particularly, topological deep learning introduced in 2017  \cite{cang2017topologynet}.

However, persistent homology has many significant drawbacks that limit its utility in applications. For example, persistent homology cannot differentiate a five-member ring from a six-member ring, such a difference is essential in molecular science. Additionally, persistent homology does not capture any non-topological changes which are crucial for networks. To overcome these drawbacks, Wei  and colleagues introduced multiscale analysis  into another fundamental object, topological Laplacians, including Hodge Laplacian and   combinatorial Laplacians in 2019, resulting in  evolutionary de Rham-Hodge theory, or persistent Hodge Laplacians  \cite{chenEvolutionaryRhamHodgeMethod2021},  and persistent combinatorial Laplacians, also called persistent  Laplacians    \cite{wang2019persistent}.  Persistent Hodge Laplacians are built from differential forms on smooth manifolds \cite{desbrun2005discrete,chenEvolutionaryRhamHodgeMethod2021}, while persistent combinatorial Laplacians are defined on simplicial complexes with point clouds \cite{horak2013spectra, chung1996combinatorial,wang2019persistent}. They are referred to as persistent topological Laplacians (PTLs) due to their similarity in algebraic topology structure. The key characteristic  of the topological Laplacian operators is that their harmonic spectra yield the Betti numbers that reveal the changes in topological invariants, while their non-harmonic spectra can convey additional information about the geometric variation \cite{horakSpectraCombinatorialLaplace2013}. By adding persistence, one can recover the persistent Betti numbers and capture the non-topological networking in the data. The development of PTLs has stimulated many new PTL formulations, including persistent sheaf Laplacians \cite{wei2021persistent}, persistent path Laplacians \cite{wang2023persistent}, and persistent hyperdigraph Laplacians \cite{chen2023persistent}, and mathematical analysis of PTLs.   M\'emoli et al. \cite{memoliPersistentLaplaciansProperties2022} developed and analyzed algorithms for computing persistent topological Laplacians, and showed a related stability theorem for some of the eigenvalues involved. Persistent Laplacians have also been generalized to simplicial maps in G\"ulen et al. \cite{gulenGeneralizationPersistentLaplacian2023}. Recently, Liu et al. have developed an interleaving distance for studying persistent Laplacians and proved some stability theorems \cite{liuAlgebraicStabilityPersistent2023}, establishing theoretical reliability of the results that persistent combinatorial Laplacians have seen in practice. Successful applications of persistent topological Laplacians have been reported in the literature, including  protein-ligand binding \cite{MengXia2021sciadv} and the accurate prediction of SARS-CoV-2 emerging dominant variants \cite{CHEN2022106262}. 
 The superiority of persistent Laplacians over classic persistent homology   was documented in a large scale machine learning-guided protein engineering study using 34 benchmark datasets \cite{qiu2023persistent}.  A software package was developed  for computing persistent Laplacians \cite{wang2021hermes}.
Recently, quantum persistent homology in terms of Dirac \cite{ameneyro2022quantum},  persistent Dirac   \cite{wee2023persistent}, and  
persistent Dirac of path and hypergraph \cite{suwayyid2023persistent}
have been proposed to achieve a similar performance for point cloud data.

Motivated by the practical need  to develop high-performance TDA models, we 
construct a directed flag Laplacian and propose persistent topological Laplacian theory on a type of simplicial complex called the directed flag complex. Directed flag complexes are based on flag complexes, also commonly referred to as clique complexes, which are an abstract simplicial complex with an underlying undirected graph structure \cite{lutgehetmannComputingPersistentHomology2020}. Flag complexes arise naturally in many problems, including computing the Vietoris-Rips complex \cite{zomorodianFastConstructionVietorisRips2010}, analyzing brain networks \cite{lutgehetmannComputingPersistentHomology2020}, and studying small world networks \cite{masulliTopologyDirectedClique2016}. Expanding the flag complex to the case of directed graphs allows researchers to encode directional information, and one can further add context from persistence that can contain geometric or time-based information. L\"ugehetmann et al. developed an open-source software FLAGSER for the fast computation of persistent homology of directed flag complexes  in 2020  \cite{lutgehetmannComputingPersistentHomology2020}. Our implementation of the present persistent directed flag Laplacian employs the underlying directed flag complex computed by FLAGSER.

The rest of this paper is organized as follows. In Section \ref{sec-background-dfl} we provide context on the directed flag complex and relevant homology tools as in \cite{lutgehetmannComputingPersistentHomology2020}. In Section \ref{sec-laplacian}, we first construct a directed flag Laplacian, and then develop a persistent directed flag Laplacian method, with some examples. The relationship of the present formalism and earlier persistent path Laplacian and persistent hyperdigraph  Laplacian are analyzed in 
Section \ref{Discussion}.  In Section \ref{sec:application} we apply the proposed persistent directed flag Laplacian to real-world data for a protein-ligand complex. A  conclusion is given in  Section \ref{sec:conclusion}.

\section{Background on Directed Flag Complex Homology}
\label{sec-background-dfl}

To establish notations,  we review the relevant graph theory and topology definitions of directed flag complexes used in \cite{lutgehetmannComputingPersistentHomology2020}. First we define the graphs we want to consider and show how to build the directed flag complex, followed by the homology of the complex. Finally, we develop persistence.

\subsection{Directed Flag Complexes}

\begin{definition}
    A directed graph, or digraph, $G$ is a pair $(V,E)$ where $V$ is a finite nonempty ordered set with elements called vertices and $E$ is a finite nonempty set of ordered pairs of vertices $e=(v_i,v_j)$ such that $v_i\neq v_j$, called edges. 
\end{definition}

Note that if $(v_i,v_j)\in E$, it is possible, but not necessarily true that $(v_j,v_i)\in E$. Also note that this definition excludes self-loops of the form $(v_i,v_i)$ from the possible edges. 

For studying the topology of these graphs, we must define additional structures on them, so we introduce the notations around abstract simplicial complexes. 

\begin{definition}
    An abstract simplicial complex on a set $X$ is a collection $\mathcal{K}$ of nonempty finite subsets of $X$, called simplices, such that if $\sigma\in\mathcal{K}$, then every nonempty $\tau\subset\sigma$ is also in $\mathcal{K}.$ An ordered simplicial complex $\mathcal{K}$ has simplices that are ordered nonempty finite subsets of $X,$ although it is not necessary that $X$ have an underlying order. If $\tau\subset \sigma$ are two simplices in a simplicial complex, we call $\tau$ a face of $\sigma$. Simplices with $k+1$ elements are called $k$-simplices.
\end{definition}

We focus our attention on a particular abstract simplicial complex of digraphs, the directed flag complex. We must specify its simplices.

\begin{definition}
    A $k$-clique of a directed graph $G=(V,E)$ is an ordered subset of vertices $\sigma=(v_0,v_1,\dots,v_{k-1})$ such that for $i<j$ we have $(v_i,v_j)\in E$. The directed flag complex ${\rm dFl}(G)$ is the ordered simplicial complex on $V$ with $k$-simplices the $(k+1)$-cliques of $G$.
\end{definition}

By this definition, all vertices of the graph are $1$-cliques and all edges are $2$-cliques. Observe that by this definition, the edges and vertices entirely determine the set of $(k+1)$-cliques, and hence the entire complex. Although every directed flag complex is an ordered simplicial complex, the requirement that all $(k+1)$-cliques are in the directed flag complex means we can construct ordered simplicial complexes, which are not directed flag complexes, by removing cliques with more than $2$ vertices.

There are many sources of directed flag complexes both in nature and in theory (the C-Elegans brain, Erd\H{o}s R\'enyi graphs) \cite{lutgehetmannComputingPersistentHomology2020}. In Example \ref{ex:example_1} we describe a few digraphs and their directed flag complexes, and in Example \ref{ex:example_2} we describe ordered simplicial complexes on the same graphs as that in Example \ref{ex:example_1}, which are not directed flag complexes.

\begin{figure}
    \centering
    \begin{subfigure}{0.3\textwidth}
        \centering
        \includegraphics[width=\textwidth]{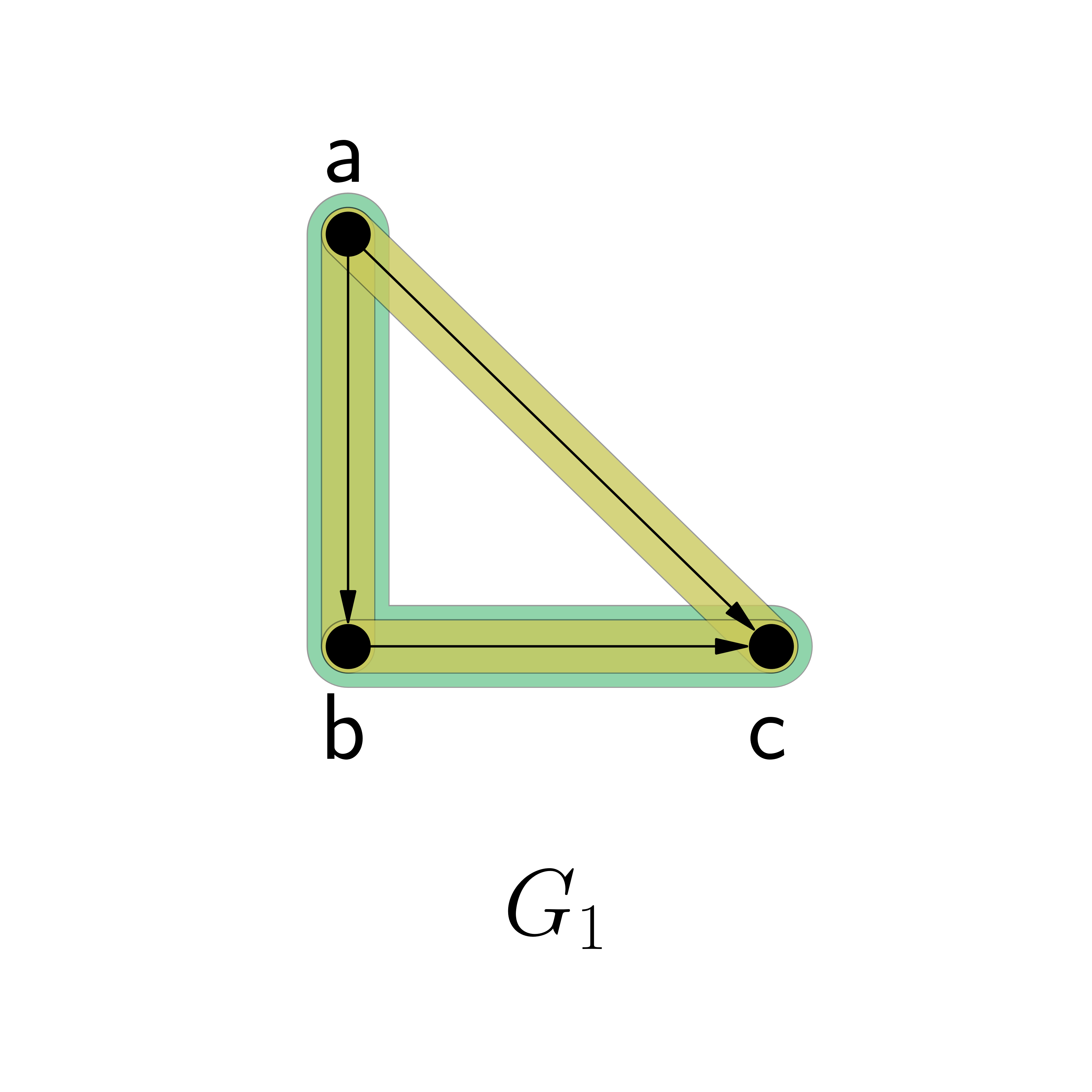}
    \end{subfigure}
    \hfill
    \begin{subfigure}{0.3\textwidth}
        \centering
        \includegraphics[width=\textwidth]{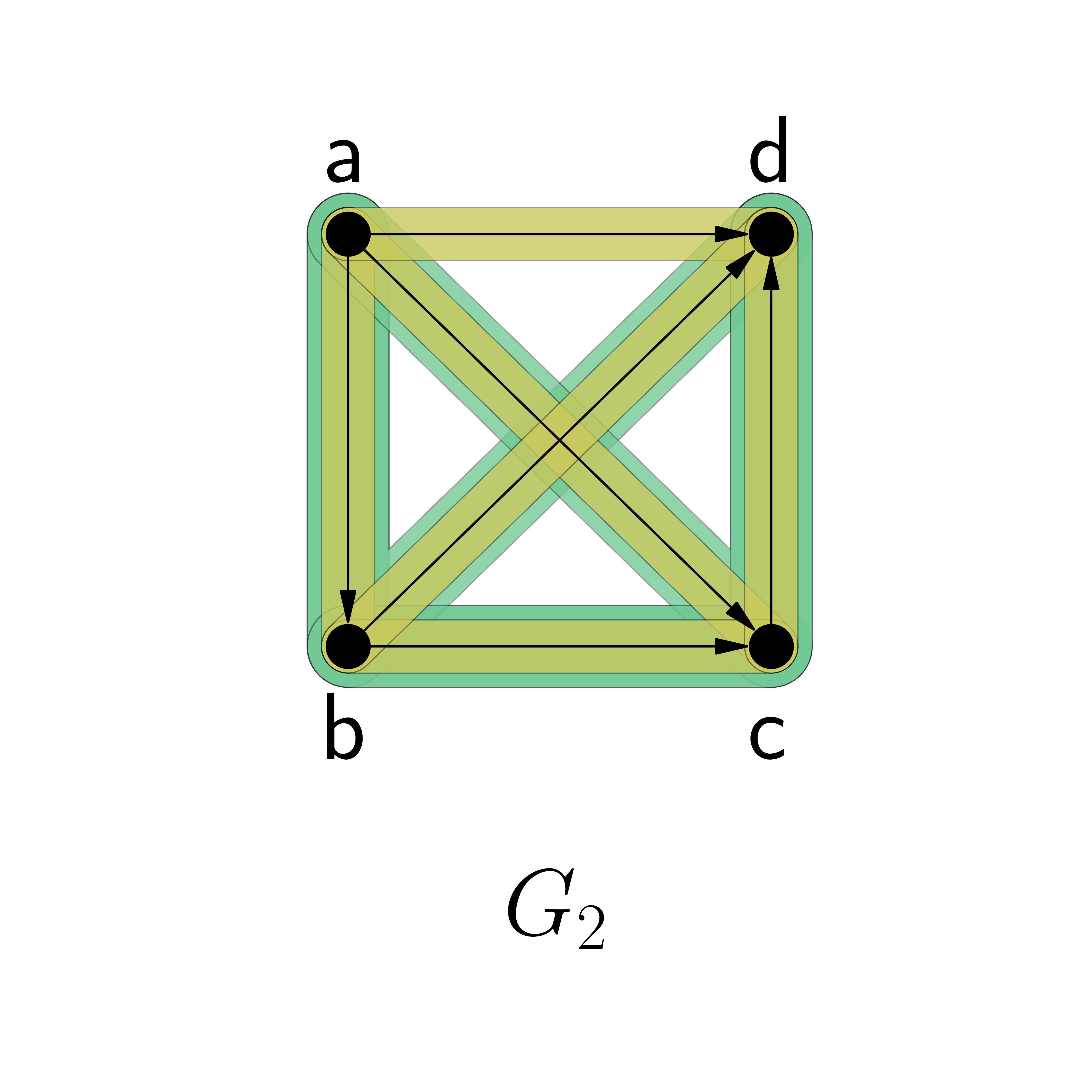}
    \end{subfigure}
    \hfill
    \begin{subfigure}{0.3\textwidth}
        \centering
        \includegraphics[width=\textwidth]{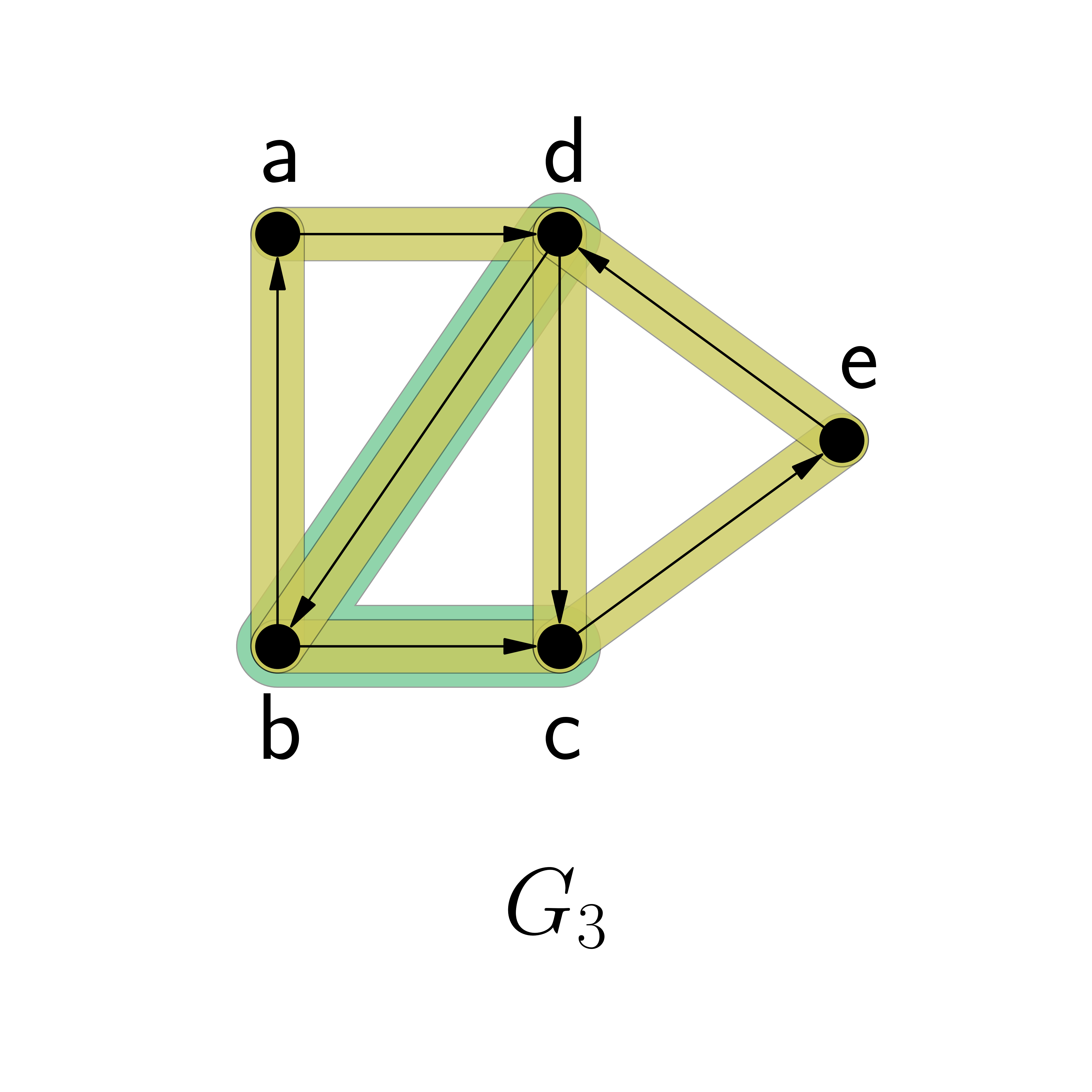}
    \end{subfigure}
    \caption{Digraphs $G_1$, $G_2,$ and $G_3.$ The edges in the directed flag complex are visualized with a yellow background around the edge, and $3$-cliques $(v_0,v_1,v_2)$ by a wider, green background around the edges $(v_0,v_1)$ and $(v_1,v_2)$, omitting $(v_0,v_2)$ to make visualization easier.}
    \label{fig:example_1}
\end{figure}

\begin{example}
\label{ex:example_1}
    Let $G_1, G_2,$ and $G_3$ be digraphs as in Figure \ref{fig:example_1}. We visualize edges in the directed flag complex with a yellow background around the edge, and $3$-cliques $(v_0,v_1,v_2)$ by a wider, green background around the edges $(v_0,v_1)$ and $(v_1,v_2)$, omitting $(v_0,v_2)$ to make visualization easier.
    
    The digraph $G_1$ has vertices (1-cliques) $V=\{a,b,c\}$, edges (2-cliques) $E=\{(a,b),(a,c),\allowbreak (b,c)\}$, and the 3-clique $(a,b,c)$. Then $dFl(G_1)$ is the ordered simplicial complex with $0$-simplices the vertices, $1$-simplices the edges, and the $2$-simplex $(a,b,c)$. 
    
    The digraph $G_2$ has $4$ vertices $V=\{a,b,c,d\}$, $6$ edges $E=\{(a,b),(a,c),(a,d),(b,c),\allowbreak (b,d),\allowbreak (c,d)\}$, $4$ $3-$cliques $\{(a,b,d),(a,b,c),(a,c,d),(b,c,d)\}$, and the $4$-clique $(a,b,c,d)$. Note that this is the minimal example of a directed $4$-clique. Then $dFl(G_2)$ is the ordered simplicial complex with $0$-simplices the vertices, $1$-simplices the edges, $2$-simplices the $4$ $3-$cliques, and $3$-simplex $(a,b,c,d)$. If we let $v_0=a, v_1=b, v_2=c, v_3=d$, then for $i<j$, $(v_i,v_j)\in E$ since $(a,b),(a,c),(a,d),(b,c),(b,d),$ and $(c,d)\in E$. 
    
    The digraph $G_3$ has $5$ vertices $\{a,b,c,d,e\}$, $7$ edges $\{(a,d),(b,a),(b,c),(c,e),(d,b),\allowbreak (d,c),\allowbreak (e,d)\}$, and $3$-clique $(d,b,c)$. Then $dFl(G_3)$ is the ordered simplicial complex with $0$-simplices the vertices, $1$-simplices the edges, and $2$-simplex $(d,b,c)$. Note that neither $(a,d,b)$ nor $(c,e,d)$ are $3$-cliques, since the edges $(a,b)$ and $(c,d)$ are not in $G_3.$
\end{example}

\begin{example}
\label{ex:example_2}
    Let $G_1, G_2,$ and $G_3$ be digraphs as in Figure \ref{fig:example_1}. Here we describe some ordered simplicial complexes that share a $1$-skeleton with the directed flag complex on these graphs, but are distinct. 

    In $G_1$, the only ordered simplicial complex with he same $1$-skeleton as the directed flag complex consists of only the vertices $\{a,b,c\}$ and edges $\{(a,b),(a,c),(b,c)\}$, with no $2$-simplices.

    In $G_2$, we can have an ordered simplicial complex with vertices $\{a,b,c,d\}$ and edges $\{(a,b),(a,c),(a,d),(b,c),(b,d),(c,d)\}$, with no $2$-simplices. We can obtain additional ordered simplicial complexes by adding any number of the possible $2$-simplices $(a,b,d),(a,b,c),\allowbreak (a,c,d)$ or $(b,c,d)$. The directed flag complex on $G_2$ is the unique ordered simplicial complex containing the $3$-simplex $(a,b,c,d)$, since all ordered subsets of $(a,b,c,d)$ must also be in the complex. By choosing $0,1,2,3,$ or all $4$ of the possible $2$-simplices on $G_2$, we can obtain $2^4 = 16$ ordered simplicial complexes that are distinct from the directed flag complex.

    If we further relax the condition that our ordered simplicial complex must share its $1$-skeleton with the directed flag complex, we can obtain many more simplicial complexes on the vertex set, although none of them are faithful to the graph. For example, on the vertices of $G_2$, we could have a simplicial complex with vertices $\{a,b,c,d\}$, edges $\{(a,b),(a,d),(d,b)\}$, and $2$-simplex $(a,d,b)$. Note that not all edges of $G_2$ must be in the complex and that the edge $(d,b)$ is not an edge of $G_2$.  
\end{example}

\subsection{Directed Flag Complex Homology}

We can define a boundary operator and chain complex on ordered simplicial complexes in the usual way.

\begin{definition}
    Let $\mathbb{F}$ be a field and $\mathcal{K}$ an ordered simplicial complex. for $k\geq 0$, let $C_k(\mathcal{K};\mathbb{F})$ be the vector space over $\mathbb{F}$ generated by $k$-simplices in $\mathcal{K}$.  The collection $C_\bullet (\mathcal{K};\mathbb{F})$ forms a chain complex with the boundary operator $d_k:C_k(\mathcal{K};\mathbb{F})\to C_{k-1}(\mathcal{K};\mathbb{F})$ given by
    \[
        d_k\left(x_0,x_1,\dots x_k\right) = \sum_{i=0}^k (-1)^i (x_0,x_1,\dots, \hat{x_i},\dots, x_k),
    \]
    where $(x_0,x_1,\dots, \hat{x_i},\dots, x_k)$ denotes the $(k-1)$-simplex obtained by removing  $x_i$ from the $k$-simplex $(x_0,x_1,\dots, x_k)$. We also take $d_0:C_0(\mathcal{K};\mathbb{F})\to 0$ to be the $0$ map. If $d_k(\sigma)=0$, we call $\sigma$ a cycle, and if $\sigma=d_{k+1}(\tau)$, we call $\sigma$ a boundary. Note that $d_k\circ d_{k-1}=0$, which also implies all boundaries are cycles.

    Since $(C_\bullet,d_\bullet)$ forms a chain complex, the $k$-th homology of $\mathcal{K}$ with coefficients in $\mathbb{F}$ is given by $H_k(\mathcal{K};\mathbb{F})=\ker(d_k)/\Ima(d_{k+1})$. Let $\beta_k=\dim H_k(\mathcal{K};\mathbb{F})$, which we call the $k$-th Betti number. 
\end{definition}

Note that in the theoretical treatment of persistent homology, the field $\mathbb{F}$ can be any field, but it is often taken to be $\mathbb{Z}/2\mathbb{Z}$ or another finite field, as is done in \cite{lutgehetmannComputingPersistentHomology2020}, for which there are fast algorithms and implementations for computing persistence \cite{edelsbrunnerTopologicalPersistenceSimplification2002}. In this work we take $\mathbb{F}$ to be the real numbers $\mathbb{R}$, because this allows us to compute the spectra of our Laplacian.

\begin{example}
    Again consider the digraph $G_3$ to be as in Figure \ref{fig:example_1}. We computed the ordered simplicial complex ${\rm dFl}(G_3)$ in Example \ref{ex:example_1}, which gives rise to the chain groups,
    
    \begin{eqnarray*}
        C_{0}({\rm dFl}(G_3);\mathbb{R}) &=& \mathrm{span}\{(a),(b),(c),(d),(e)\}, \\
        C_{1}({\rm dFl}(G_3);\mathbb{R}) &=& \mathrm{span}\{(a,d), (b,a), (b,c), (c,e), (d,b), (d,c), (e,d)\}, \\
        C_{2}({\rm dFl}(G_3);\mathbb{R}) &=& \mathrm{span}\{(d,b,c)\}.
    \end{eqnarray*}
    
    The nonzero boundary operators can then be represented by the matrices
    
    \begin{align*}
        d_1 &= \begin{blockarray}{cccccccc}
            & (a,d) & (b,a) & (b,c) & (c,e) & (d,b) & (d,c) & (e,d) \\
            \begin{block}{c(ccccccc)}
              (a) & -1 & 1  & 0  & 0  & 0  & 0 &0 \\
              (b) & 0  & -1 & -1 & 0  & 1  & 0 & 0 \\
              (c) & 0  & 0  & 1  & -1 & 0  & 1 & 0 \\
              (d) & 1  & 0  & 0  & 0  & -1 & -1 & 1 \\
              (e) & 0  & 0  & 0  & 1  &  0 & 0 & -1 \\
        \end{block}
        \end{blockarray}\\  
        d_2 &= \begin{blockarray}{cc}
            & (d,b,c)\\
            \begin{block}{c(c)}
                (a,d) & 0\\
                (b,a) & 0\\
                (b,c) & 1\\
                (c,e) & 0\\
                (d,b) & 1\\
                (d,c) & -1\\
                (e,d) & 0\\
            \end{block}
        \end{blockarray}.
    \end{align*}

    We can then compute the kernels and images of the relevant boundary maps,

    \begin{eqnarray*}
        \ker d_0 &=& \mathrm{span}\{(a),(b),(c),(d),(e)\}, \\
        \ker d_1 &=& \mathrm{span}\{(a,d)+(d,b)+(b,a),(d,b)-(d,c)+(b,c),(c,e)+(e,d)+(d,c)\},\\
        \Ima d_1 &=& \mathrm{span}\{(a)-(b),(c)-(b),(d)-(a),(e)-(c)\},\\
        \Ima d_2 &=& \mathrm{span}\{(b,c)-(d,c)+(d,b)\}
    \end{eqnarray*}.

    Note that $\Ima d_0, \ker d_2,\ker d_i,$ and $\Ima d_i$ are all $0$ for $i\geq 3$.
    From these, we compute the homology,
    \begin{align*}
        H_0({\rm dFl}(G);\mathbb{R}) = \ker d_0/\Ima d_1 &\cong \mathbb{R},\\
        H_1({\rm dFl}(G);\mathbb{R}) = \ker d_1/\Ima d_2 &\cong \mathbb{R}^2,\\
        H_2({\rm dFl}(G);\mathbb{R}) = \ker d_2/\Ima d_3 &\cong 0.
    \end{align*}

    All other homologies are $0$. Then the Betti numbers are $\beta_0=1, \beta_1 = 2,$ and $\beta_2=0$.
\end{example} 

This example demonstrates a directed flag complex with a cycle in $C_1$ that is a boundary and two that are not boundaries. It also demonstrates that as usual, $\beta_0$ measure the number of connected components.

\subsection{Persistent Directed Flag Complex Homology}\label{subsec:persistent_flag_homology}

Now that we have established homology of directed flag complexes, we introduce the relevant details of persistence, including the persistent Betti numbers. Let $(\mathbb{R},\leq)$ be the category of real numbers with morphisms given by $a\to b$ if $a\leq b$. Let $\dflcat$ be the category of directed flag complexes, with morphisms given by inclusion. A filtration is a functor $\mathcal{F}:(\mathbb{R},\leq)\to \dflcat$, where $\mathcal{F}(a\leq b)$ is the inclusion $i^{a,b}:\mathcal{F}(a)\hookrightarrow \mathcal{F}(b).$ For clarity of notation, we denote $H_k^a(\mathcal{F})=H_k(\mathcal{F}(a),\mathbb{R})$ and the boundary map of this chain complex to be $d_k^a:H_k^a\to H_{k-1}^a$. When the filtration $\mathcal{F}$ is understood from context, we use $H_k^a$. Then for a given filtration $\mathcal{F}$ and dimension $k$, we define the real $(a,b)-$persistent homology to be

\[
    H_k^{a,b}(\mathcal{F}) = \Ima\left(H_k^a\to H_k^b\right),
\]

and the $(a,b)$-persistent Betti number of $\mathcal{F}$ to be $\beta_k^{a,b}=\dim H_k^{a,b}(\mathcal{F})$. Again, if the filtration is understood, we use $H_k^{a,b}$. We note that the coefficients do not need to be $\mathbb{R}$ in general, but for the development of the persistent Laplacian we will use $\mathbb{R}.$ 

We similarly denote $C_k^a(\mathcal{F})=C_k(\mathcal{F}(a);\mathbb{R})$ or simply $C_k^a$ if the filtration is clear from the context. The persistent Betti number $\beta_k^{a,b}$ corresponds to the number of nonzero classes in $H_k^a$ that are nonzero in $H_k^b$.

We now give an example where persistent directed flag complex homology is different from computing the non-persistent directed flag complex homology on a sequence of graphs.

\begin{figure}[!ht]
    \centering
    \includegraphics[width=0.5\textwidth]{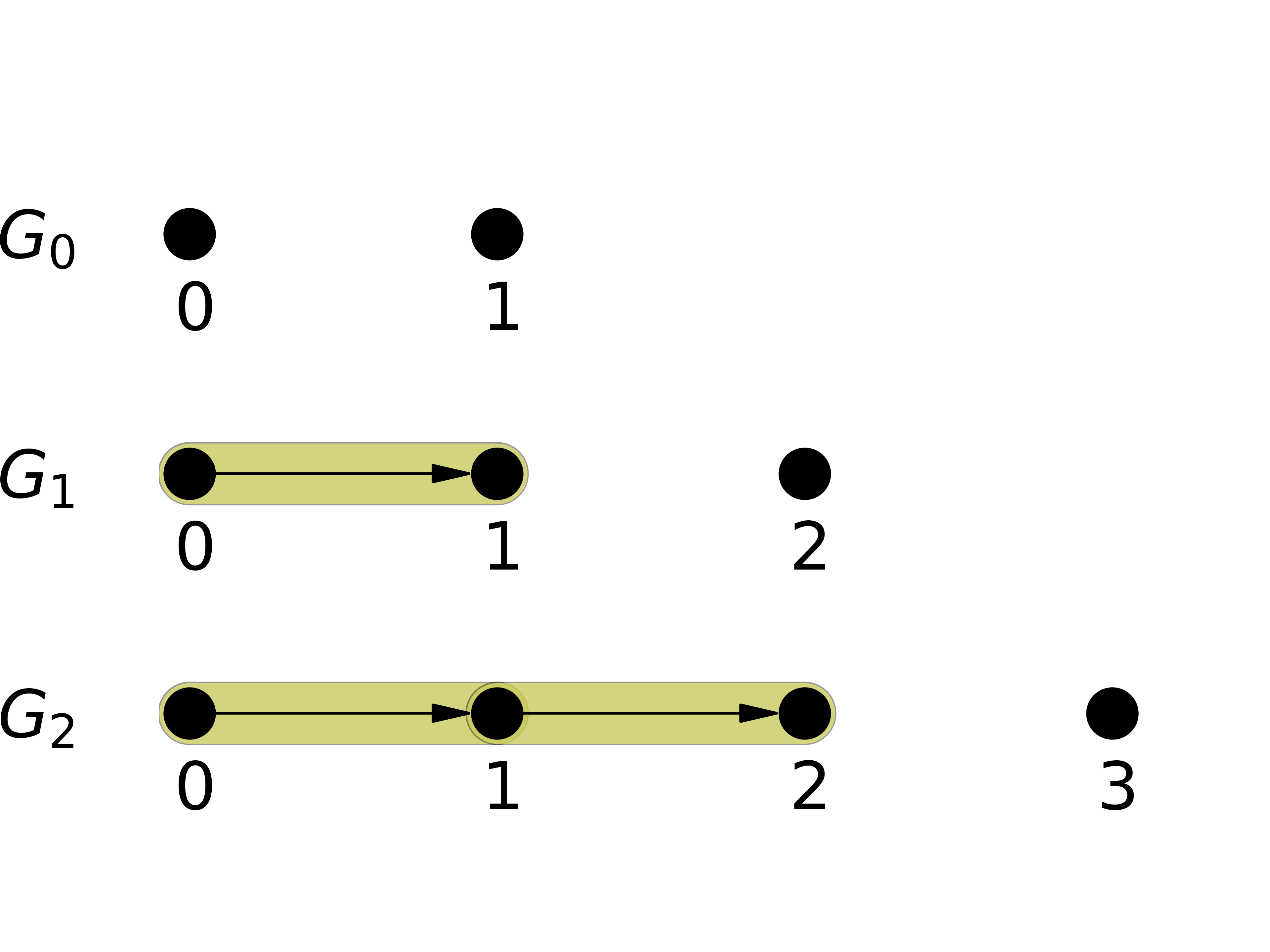}
    \caption{Filtered directed flag complex with constant persistent Betti number $\beta_0^{a,a+1}=1$, while the non-persistent Betti number is constant at $\beta_0^a=2$. The area around edges are colored in yellow to signify that they form $1$-simplices. }
    \label{fig:const_betti_0}
\end{figure}

\begin{example}
In Fig. \ref{fig:const_betti_0} we show a filtration of graphs, where the persistent directed flag complex homology has $\beta_0^{0,1}=\beta_0^{1,2}=1$, while $\beta_0^0=\beta_0^1=\beta_0^2=2$.  This can be seen by the fact that only $1$ connected component that existed in $G_0$ also exists in $G_1$, hence $\beta_0^{0,1}=1$, despite there being $2$ total connected components in $G_1$.

Alternatively, consider the induced map by inclusion $i_0^{1,2}:H_0^1\to H_0^2$, where $H_0^1\cong\mathrm{span}\{(0),(1)\}$ and $H_0^2\cong\mathrm{span}\{(0),(2)\}$. Then $\Ima i_0^{1,2}\cong\mathrm{span}\{(0)\}\cong\mathbb{R}$, so $\beta_0^{1,2} = 1$. Similar computations apply to $\beta_0^{2,3}$.
\end{example}

This example demonstrated that persistent Betti numbers provide a different perspective than calculating a sequence of non-persistent Betti numbers. Moreover, since the non-persistent Laplacian spectra recover the non-persistent Betti numbers, this motivates the need for the persistent Laplacian to recover the persistent Betti numbers.  

\section{Persistent Directed Flag Laplacian}
\label{sec-laplacian}

Persistent Betti numbers do not tell the full story of the persistence evolution. Another way to measure the changes with persistence is through persistent topological Laplacians like the persistent path Laplacian \cite{wang2023persistent}, topological hyperdigraph Laplacian\cite{chen2023persistent}, and persistent sheaf Laplacian \cite{wei2021persistent}. In this work, we extend the persistent directed flag complex homology to a persistent directed flag Laplacian (PDFL).  
To do this, we construct the directed flag Laplacian and calculate its spectra. Then we construct the persistent version of the directed flag Laplacian and calculate its spectra.

\subsection{Directed Flag Laplacian}
\label{subsec-dfl}
In this section, let $G=(V,E)$ be a digraph, ${\rm dFl}(G)$ be the directed flag complex of $G$, and $S_k$ be the set of $k$-cliques in $G$. Also write $C_k=C_k({\rm dFl}(G);\mathbb{R}).$ The Laplacian requires an adjoint of the boundary operator, so we define an inner product on the basis $S_k$ of $C_k$,
\[
    \langle\cdot,\cdot\rangle:S_n\times S_n\to \mathbb{R},
\]

by 

\[
    \langle\sigma,\tau\rangle= \begin{cases}1 & \text{if }\sigma=\tau,\\
                                            0 & \text{otherwise}
                                \end{cases}.
\]

Next, this induces an inner product $\llangle\cdot,\cdot\rrangle$ on the dual space $C^k=C^k({\rm dFl}(G);\mathbb{R})=\Hom(C_k,\mathbb{R})$ such that for any $f,g\in C^k$,

\[
    \llangle f,g\rrangle = \sum_{\sigma \in S_k} f(\sigma)g(\sigma).
\]

With respect to these inner products, there is an adjoint $(d_k)^*:C_{k-1}\to C_k$ of $d_k$. Then we may define the $k$-th combinatorial directed flag Laplacian, or directed Flag Laplacian, by $\Delta_k:C_k\to C_k$,

\[
    \Delta_k = d_{k+1}\circ \left(d_{k+1}\right)^* + \left(d_k\right)^* \circ d_k,
\]

where $\Delta_0= d_1 d_1^*$ and if $N$ is the largest integer for which $S_N$ is nonempty, $\Delta_N = d_N^* d_N$. With respect to the standard basis of $C_k$, let $B_k$ be the matrix representation of $d_k$. In the unweighted, non-persistent case, as shown in \cite{memoliPersistentLaplaciansProperties2022}, the matrix representation of $(d_k)^*$ with respect to the standard basis of $C_{k-1}$ is $B_k^T$. In the computing the matrix representation of the adjoint in the persistent case, we will need more care as shown in \cite{memoliPersistentLaplaciansProperties2022}. The matrix representation of the $k$-th combinatorial directed flag Laplacian $\Delta_k:C_k\to C_k$ is 

\[
    L_k = B_{k+1}B_{k+1}^T + B_k^T B_k.
\]

Then $L_k$ is positive semi-definite, with only real, non-negative eigenvalues, and is self-adjoint. Moreover, $\beta_k=\dim\ker\Delta_k$, which is the number of zero eigenvalues of $L_k$, as shown in \cite{wang2019persistent}. 

We may extract more information from $L_k$ than just the multiplicity of zero as an eigenvalue by examining the nonzero eigenvalues, called the non-harmonic spectra. If the dimension of $C_k$ is $n$, denote the spectra of $\Delta_k$ in non-decreasing order by

\[
    \mathrm{Spectra}(L_k) = \{(\lambda_1)_k,(\lambda_2)_k,\dots,(\lambda_n)_k\}.
\]

We now compute the Laplacians and spectra for the directed flag complex in Example \ref{ex:example_1}, and show that the Betti numbers we previously computed agree with the harmonic spectra of the Laplacians.

\begin{example}
    Using the graph and directed flag complex  of $G_3$ in Example \ref{ex:example_1}, we compute the Laplacians and spectra. The nonzero boundary maps from before are,
    
	 \begin{align*}
        B_1 &= \begin{blockarray}{cccccccc}
            & (a,d) & (b,a) & (b,c) & (c,e) & (d,b) & (d,c) & (e,d) \\
            \begin{block}{c(ccccccc)}
              (a) & -1 & 1  & 0  & 0  & 0  & 0 &0 \\
              (b) & 0  & -1 & -1 & 0  & 1  & 0 & 0 \\
              (c) & 0  & 0  & 1  & -1 & 0  & 1 & 0 \\
              (d) & 1  & 0  & 0  & 0  & -1 & -1 & 1 \\
              (e) & 0  & 0  & 0  & 1  &  0 & 0 & -1 \\
        \end{block}
        \end{blockarray}\\  
        B_2 &= \begin{blockarray}{cc}
            & (d,b,c)\\
            \begin{block}{c(c)}
                (a,d) & 0\\
                (b,a) & 0\\
                (b,c) & 1\\
                (c,e) & 0\\
                (d,b) & 1\\
                (d,c) & -1\\
                (e,d) & 0\\
            \end{block}
        \end{blockarray}.
    \end{align*}

    Then

    \begin{align*}
        L_0 &= \begin{pmatrix}
			-1 & 1  & 0  & 0  & 0  & 0 & 0 \\
			0  & -1 & -1 & 0  & 1  & 0 & 0 \\
			0  & 0  & 1  & -1 & 0  & 1 & 0 \\
			1  & 0  & 0  & 0  & -1 & -1 & 1 \\
			0  & 0  & 0  & 1  &  0 & 0 & -1 
        \end{pmatrix}
        \begin{pmatrix}
            -1 & 0  & 0  & 1  & 0 \\
			1  & -1 & 0  & 0  & 0 \\
			0  & -1 & 1  & 0  & 0 \\
			0  & 0  & -1 & 0  & 1 \\
			0  & 1  & 0  & -1 & 0\\
			0 &  0 &  1 & -1 &  0\\
 			0 &  0 &  0 &  1 & -1
        \end{pmatrix}\\
        &= \begin{pmatrix} 
             2 &-1 & 0 & -1 & 0\\
            -1 & 3 &-1 &-1 & 0\\
             0  & -1 & 3 & -1 &-1\\
             -1 &-1 &-1 & 4 &-1\\
             0 & 0 &-1 &-1 & 2
        \end{pmatrix},\\
        L_1 &= \begin{pmatrix}
            0\\
            0\\
            1\\
            0\\
            1\\
            -1\\
            0
        \end{pmatrix}
        \begin{pmatrix}
            0 & 0 & 1 & 0 & 1 & -1 & 0
        \end{pmatrix}
        +
        \begin{pmatrix}
            -1 & 0  & 0  & 1  & 0 \\
			1  & -1 & 0  & 0  & 0 \\
			0  & -1 & 1  & 0  & 0 \\
			0  & 0  & -1 & 0  & 1 \\
			0  & 1  & 0  & -1 & 0\\
			0 &  0 &  1 & -1 &  0\\
 			0 &  0 &  0 &  1 & -1
        \end{pmatrix}\begin{pmatrix}
			-1 & 1  & 0  & 0  & 0  & 0 & 0 \\
			0  & -1 & -1 & 0  & 1  & 0 & 0 \\
			0  & 0  & 1  & -1 & 0  & 1 & 0 \\
			1  & 0  & 0  & 0  & -1 & -1 & 1 \\
			0  & 0  & 0  & 1  &  0 & 0 & -1 
        \end{pmatrix}\\
        &= \begin{pmatrix} 
             2 &-1 & 0 &0 &-1 & -1 & 1\\
            -1 & 2 & 1 &0 & -1 & 0 & 0\\
             0 & 1 & 3 & -1 & 0 & 0 & 0\\
            0 & 0 & -1 & 2 & 0 & -1 &-1\\
            -1 & -1 & 0 & 0 & 3 & 0 &-1\\
             -1 & 0 & 0 & -1 & 0 & 3 & -1\\
             1 & 0 & 0 & -1 & -1 & -1 & 2
        \end{pmatrix},\\
        L_2 &= \begin{pmatrix}
            0 & 0 & 1 & 0 & 1 & -1 & 0
        \end{pmatrix}
        \begin{pmatrix}
            0\\
            0\\
            1\\
            0\\
            1\\
            -1\\
            0
        \end{pmatrix}\\
        &= \begin{pmatrix}
            3
        \end{pmatrix}.   
    \end{align*}

    From these matrices we can compute the spectra:

    \begin{align*}
        \mathrm{Spectra}(L_0) &= \{0,3-\sqrt{2},3, 3+\sqrt{2},5\}\\
        \mathrm{Spectra}(L_1) &= \{0,0,3-\sqrt{2},3,3,3+\sqrt{2},5\}\\
        \mathrm{Spectra}(L_2) &= \{3\}.
    \end{align*}

    Note that some eigenvalues are repeated. From these spectra, we can see that $\beta_0=1,\beta_1=2,$ and $\beta_2=0$, which agrees with our prior calculations.
    
\end{example}

\subsection{Persistent Directed Flag Laplacian}

We now move to the persistent case, where we have a filtration functor $\mathcal{F}:(\mathbb{R},\leq)\to\dflcat$. To do this, we repeat the process of construction as in Section \ref{subsec-dfl} with the added structure of a filtration. For building the persistent directed flag Laplacian, we will work more directly with the chain groups than is necessary in formulating persistent homology as in Section \ref{subsec:persistent_flag_homology}. To work with the chain groups, we note that for $k\geq 0$, the inclusion of directed flag complexes induces morphisms on chains $i^{a,b}_k:C_k^a\to C_k^b$, which also induces a morphism on homologies, $i_k^{a,b}:H_k^a\to H_k^b$. 

Let $C_{k+1}^{a,b}=\{\sigma\in C_{k+1}^b\mid d_{k+1}^b \sigma \in C_{k}^a\subset C_{k}^b\}$ and $\iota_{k+1}^{a,b}$ the inclusion $C_{k+1}^{a,b}\hookrightarrow C_{k+1}^b$. Then $C_{k+1}^{a,b}$ can be described as the chains in $C_{k+1}^{b}$ with boundaries contained entirely in $C_k^a$. It would otherwise be possible that some element of $C_{k+1}^b$ has boundary intersecting $C_{k}^b\setminus C_k^a$. Then we define the persistent boundary operator $d_{k+1}^{a,b}=(i_k^{a,b})^*\circ d_{k+1}^b\circ\iota_{k+1}^{a,b}$. The following commutative diagram shows how the relevant maps fit into the context of the chain complexes.

\begin{center}
\begin{tikzcd}
C_{k+1}^a \arrow[dd, "{i_{k+1}^{a,b}}"', hook, shift left] \arrow[rr, "d_{k+1}^a"] &                                                                                                      & C_k^a \arrow[rr, "d_{k}^a", shift left] \arrow[ld, "{\left(d_{k+1}^{a,b}\right)^*}"] \arrow[dd, "{i_k^{a,b}}", hook, shift left] &  & C_{k-1}^a \arrow[ll, "\left(d_{k+1}^a\right)^*", shift left] \arrow[dd, "{i_{k-1}^{a,b}}", hook, shift left] \\
                                                                                   & {C_{k+1}^{a,b}} \arrow[ld, "{\iota_{k+1}^{a,b}}"', hook] \arrow[ru, "{d_{k+1}^{a,b}}", shift left=2] &                                                                                                                                  &  &                                                                                                              \\
C_{k+1}^b \arrow[rr, "d_{k+1}^b"]                                                  &                                                                                                      & C_k^b \arrow[rr, "d_k^b"]                                                                                                        &  & C_{k-1}^b                                                                                                   
\end{tikzcd}
\end{center}

Define the $k$-th $(a,b)$-persistent directed flag Laplacian $\Delta_k^{a,b}:C_k^a\to C_k^b$ by
\[
    \Delta_k^{a,b} = d_{k+1}^{a,b}\circ \left(d_{k+1}^{a,b}\right)^* + \left(d_k^a\right)^*\circ d_k^a.
\]

In the literature, the terms $d_{k+1}^{a,b}\circ \left(d_{k+1}^{a,b}\right)^*$ and $\left(d_k^a\right)^*\circ d_k^a$ are sometimes called the ``up" and ``down" Laplacians, respectively. A detailed proof that the harmonic spectra of $\Delta_k^{a,b}$ correspond to the $\beta_k^{a,b}$ can be found in Theorem 2.7 of \cite{memoliPersistentLaplaciansProperties2022}.

We need a matrix representation $B_{k+1}^{a,b}$ of $d_{k+1}^{a,b}$ with respect to a basis of $C_{k+1}^{a,b}$. The choice of basis matters, as discussed in \cite{memoliPersistentLaplaciansProperties2022}, because the matrix representation of $\left(d_{k+1}^{a,b}\right)^\ast$ depends on the inner product structure of $C_{k+1}^{a,b}$, which is an inner product subspace of $C_{k+1}^b$, and $C_{k+1}^{a,b}$ does not necessarily have a basis that is a restriction of the standard basis of $C_{k+1}^{a,b}$.

If $Z$ is a matrix with columns that form a basis of $C_{k+1}^{a,b}$ in $C_{k+1}^b$ and $B_{k+1}^{a,b}$ is written with respect to that basis, then the matrix representation of $\left(d_{k+1}^{a,b}\right)^\ast$ is $(Z^TZ)^{-1}\left(B_{k+1}^{a,b}\right)^T$. If we take our basis to be orthonormal, then $Z^TZ$ is the identity matrix on $C_{k+1}^{a,b}$, so the matrix representation of  $\left(d_{k+1}^{a,b}\right)^\ast$ is simply the transpose $\left(B_{k+1}^{a,b}\right)^T$. Algorithm 3.1 of \cite{memoliPersistentLaplaciansProperties2022} provides an explicit procedure that can be used for finding $Z$. We slightly modify this procedure by transforming the basis $Z$ into an orthonormal one. This gives a matrix representation of $\iota_{k+1}^{a,b}$ with respect to an orthonormal basis of $C_{k+1}^{a,b}$ and the standard basis of $C_{k+1}^b$. The matrix representation $J_k^{a,b}$ of $i_{k}^{a,b}$ is taken to be the identity on elements of $C_k^a$ and $0$ otherwise. 

This gives us a matrix representation $B_{k+1}^{a,b}=J_{k}^{a,b}B_{k+1}^bZ$ of $d_{+1}^{a,b}$ and a representation $\left(B_{k+1}^{a,b}\right)^T$ of the adjoint $\left(d_{+1}^{a,b}\right)^\ast$. Algorithm 3.1 of \cite{memoliPersistentLaplaciansProperties2022} uses a series of row and column operations that are equivalent to the matrix multiplication by $Z$ and $J_k^{a,b}$.  Replacing the computation of $Z^TZ$ in algorithm 3.1 of \cite{memoliPersistentLaplaciansProperties2022} with computing an orthonormal basis of $Z$ is an asymptotically equal tradeoff in complexity. In practice, the slowest steps of computing the persistent directed flag Laplacian seem to be in building the complex and performing the column reduction over $\mathbb{R}$ for finding $B_{k+1}^{a,b}$. 

The current work focuses primarily on the knowledge gained from extending persistent direct flag homology to a Laplacian theory rather than the algorithmic efficiency. Because our goal is theoretical in nature, multiple avenues for improvement in speed have been left on the table. We note that Algorithm 4.1 of \cite{memoliPersistentLaplaciansProperties2022} provides an alternative method to compute the persistent Laplacian using Schur complements, which the authors explain may have a better time complexity in special cases, including that of clique (flag) complexes.

Then the matrix representation of $\Delta_k^{a,b}$ is
\[
    L_k^{a,b} = \left(B_{k+1}^{a,b}\right)^T B_{k+1}^{a,b} + B_k^a \left(B_k^a\right)^T
\]

and we denote its spectra in non-decreasing order as

\[
    \mathrm{Spectra}(L_k^{a,b}) = \{(\lambda_1)_k^{a,b},(\lambda_2)_k^{a,b},\dots,(\lambda_n)_k^{a,b}\}.
\]

The smallest nonzero eigenvalue of $L_k^{a,b}$ is often of importance, so we denote it $\lambda_k^{a,b}$.

\begin{figure}[!ht]
    \centering
    \includegraphics[width=0.6\textwidth]{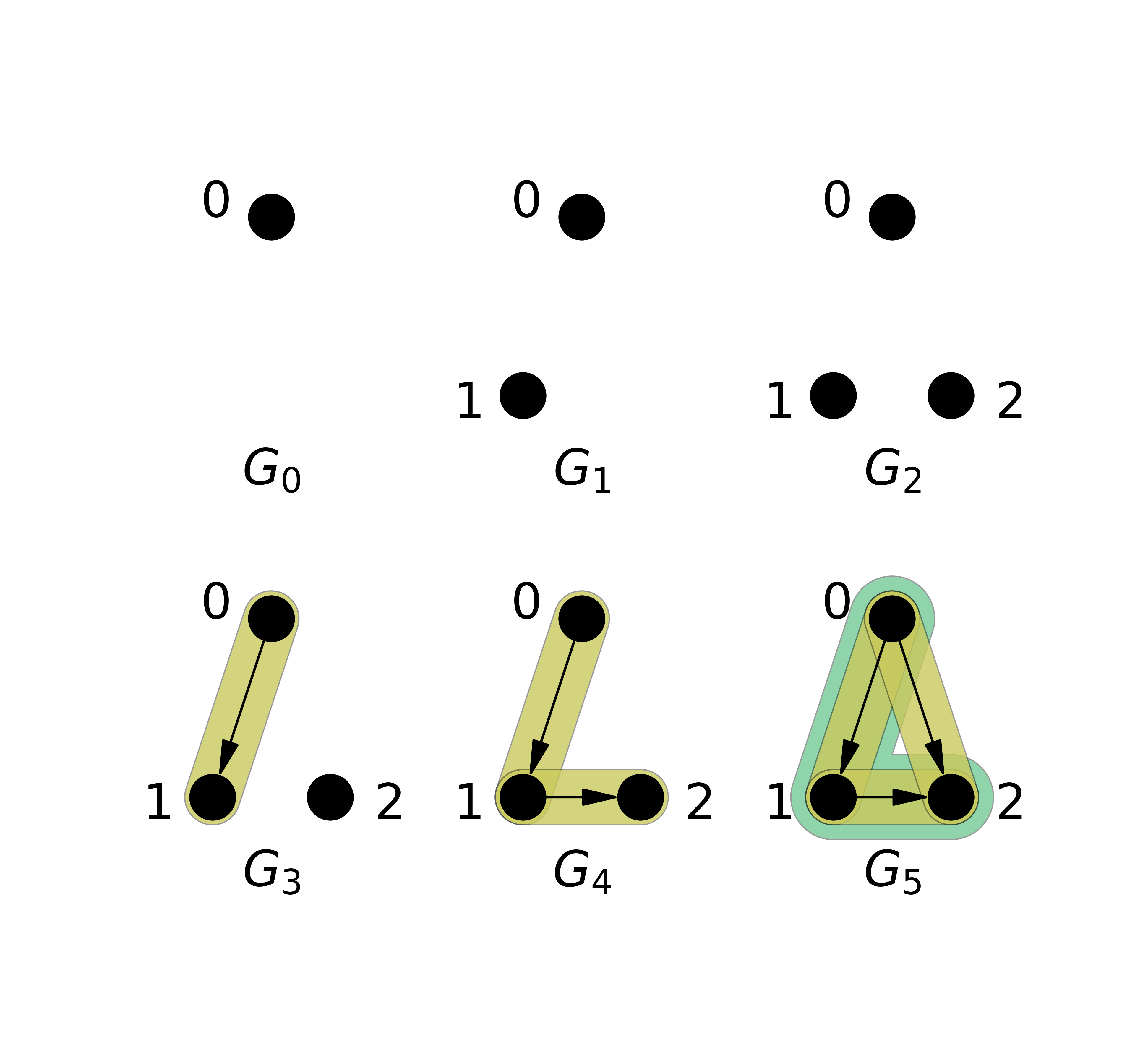}
        \caption{Filtration of a directed flag complex on a graph with three vertices. The area around edges are colored in yellow to signify that they form $1$-simplices, and the area around the edges $(v_0,v_1)$ and $(v_1,v_2)$ is shaded green to highlight them as part of the $2$-simplex $(v_0,v_1,v_2)$.}
    \label{fig:filtered_triangle}
\end{figure}

\begin{figure}[!ht]
    \centering
    \includegraphics[width=0.8\textwidth]{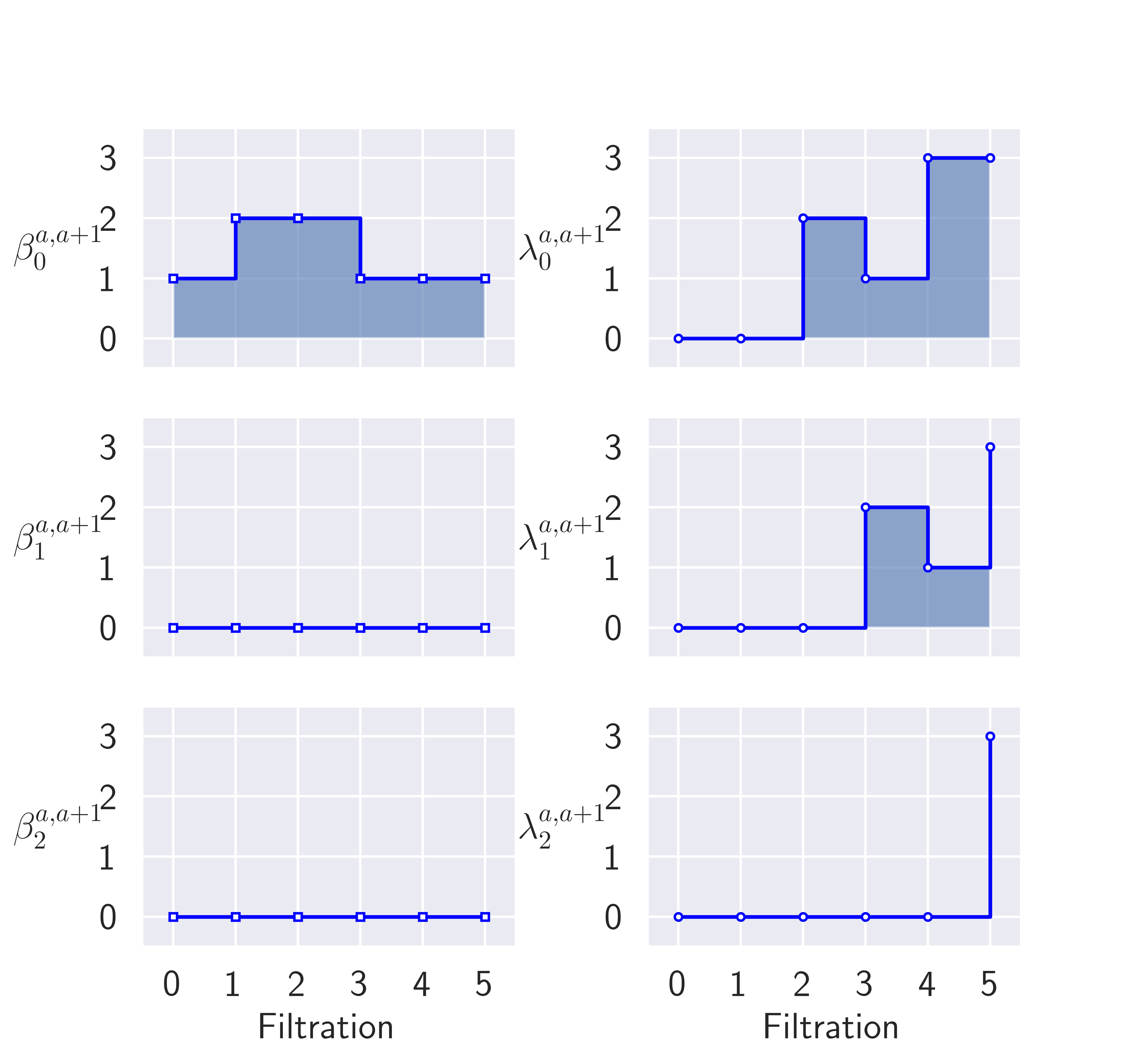}
    \caption{Persistent spectra of the graph in Fig. \ref{fig:filtered_triangle}. The values of $\beta_i^{5,-}$ and $\lambda_i^{5,-}$ are $\beta_i^5$ and $\lambda_i^5$, respectively, as filtered complex does not change after this point. }
    \label{fig:filtered_triangle_spectra}
\end{figure}

In Fig. \ref{fig:filtered_triangle}, we display a filtration of directed flag complexes on a graph with three vertices, and in Fig. \ref{fig:filtered_triangle_spectra}, we show the persistent Betti numbers and $\lambda_k^{a,a+1}$. When the edge $(0,2)$ is added in $G_5$, the directed flag complex also gains a $2$-cell, $(0,1,2).$ There is a cycle created in $C_1^5$, $(1,2)-(0,2)+(0,1)$, but it is also the boundary of $(0,1,2)\in C_2^5$. This cycle does not generate a homology class in the first homology group, so persistent Betti numbers cannot detect it, but $\lambda_1^{a,a+1}$ shows a change.

\section{Discussion}\label{Discussion}

Here we place the directed flag Laplacian we have built in context with the path Laplacian described in \cite{wang2023persistent}, based on path homology as described in \cite{GRIGORYAN2019106877}, and with the hyperdigraph Laplacian described in \cite{chen2023persistent}. These other two topological Laplacians both generalize simplicial complexes, and hence generalize the directed flag complex, but in interesting ways.

\subsection{Relationship with the path Laplacian}

Here we include some of the key components of the path Laplacian theory. The full details can be found in \cite{wang2023persistent}. First, we define the terminology around paths. For a nonempty finite set $V$, with elements called vertices, and nonnegative integer $p$, an elementary $p-$path on $V$ is any sequence $i_0\dots i_p$ of $p+1$ vertices. Let $\Lambda_p$ be the vector space over a field $\mathbb{K}$ generated by the elementary $p$-paths, with elements called $p$-paths and the elements representing the elementary $p$-paths are denoted by $e_{i_0\dots i_p}$. The authors define the elementary $(-1)-$path to be the empty set, so the $\Lambda_{-1}\cong\mathbb{K}$, and also define $\Lambda_{-2}=\{0\}$. A boundary operator $\partial$ defined in the usual way on $\Lambda_{\ast}$ produces a chain complex. 

A path complex on $V$ is a nonempty collection $P$ of elementary paths such that 

\begin{equation}\label{eq:path_complex}
    \text{if } i_0\dots i_n\in P,\text{ then }i_0\dots i_{n-1}\in P\text{ and }i_1\dots i_{n}\in P.
\end{equation}

The set of $n$-paths in $P$ is denoted $P_n$. 

 The allowed $n$-paths are given by
\[
\mathcal{A}_n=\mathrm{span}_\mathbb{K}\{e_{i_0\dots i_n} : i_0\dots i_n\in P_n\},
\]
which is a subspace of $\Lambda_n$. In order to form a good chain complex of allowed paths, one must further restrict to $\partial$-invariant paths, 
\[
    \Omega_n=\{v\in \mathcal{A}_n : \partial v\in \mathcal{A}_{n-1}\}.
\]

The path homology groups $H_n(P)$ of the path complex $P$ are given by the homology of the chain complex

\[
    \Omega_n \xrightarrow[]{\partial}\Omega_{n-1}\xrightarrow[]{\partial}\cdots\xrightarrow[]{\partial}\Omega_0\xrightarrow[]{\partial}0.
\]

By taking $\mathbb{K}=\mathbb{R}$ and defining an appropriate inner product, one can obtain the adjoint $\partial_n^\ast$ and construct a Laplacian theory for path complexes in a manner similar to what we have done with the directed flag complex. Denote the path Laplacian by $\Delta_{n}^{P}=\partial_{n+1}^{P}\left(\partial_{n+1}^{P}\right)^{\ast} + \left(\partial_{n}^{P}\right)^{\ast}\partial_{n}^{P}$, with corresponding matrix representation $L_n^{P}=B_{n+1}^P \left(B_{n+1}^P\right)^T + \left(B_n^P\right)^T B_n^p$. It is pointed out in \cite{memoliPersistentLaplaciansProperties2022} that the transpose is not always the correct matrix representation for the adjoint of the persistent boundary matrix; the same is true here, even in the non-persistent case, because $\Omega_n$ is an inner product subspace of $\mathcal{A}_n$ with a basis which may not always be written as a subset of the basis for $\mathcal{A}_n$. As with the directed flag Laplacian, we will write $B_{n}^{P}$ with respect to an orthonormal basis for $\Omega_n$, which resolves this issue.

It is noted in \cite{wang2023persistent} and \cite{grigoryanPathComplexesTheir2020} that path complexes are designed to generalize abstract simplicial complexes, and hence any directed flag complex can be given an equivalent path complex. Note that $P_n$ does not necessarily contain all elementary $n$-paths on $V$; it only needs to satisfy the conditions in Eq. \ref{eq:path_complex}. Although path complexes generalize abstract simplicial complexes, and hence the directed flag complex, the most natural path complex to use on a graph $G$ will not necessarily agree with the directed flag complex on $G$.

For a digraph $G=(V,E)$, it is most natural to consider the path complex $P(G)$ consisting of all vertices, edges, and paths that can be built with those edges. Additional restrictions on the digraph are discussed in \cite{wang2023persistent}, but the present level of detail suffices to distinguish the directed flag Laplacian from the path Laplacian. This path complex will generally not agree with the directed flag complex on the same graph.

\begin{figure}[!ht]
    \centering
    \includegraphics[width=0.6\textwidth]{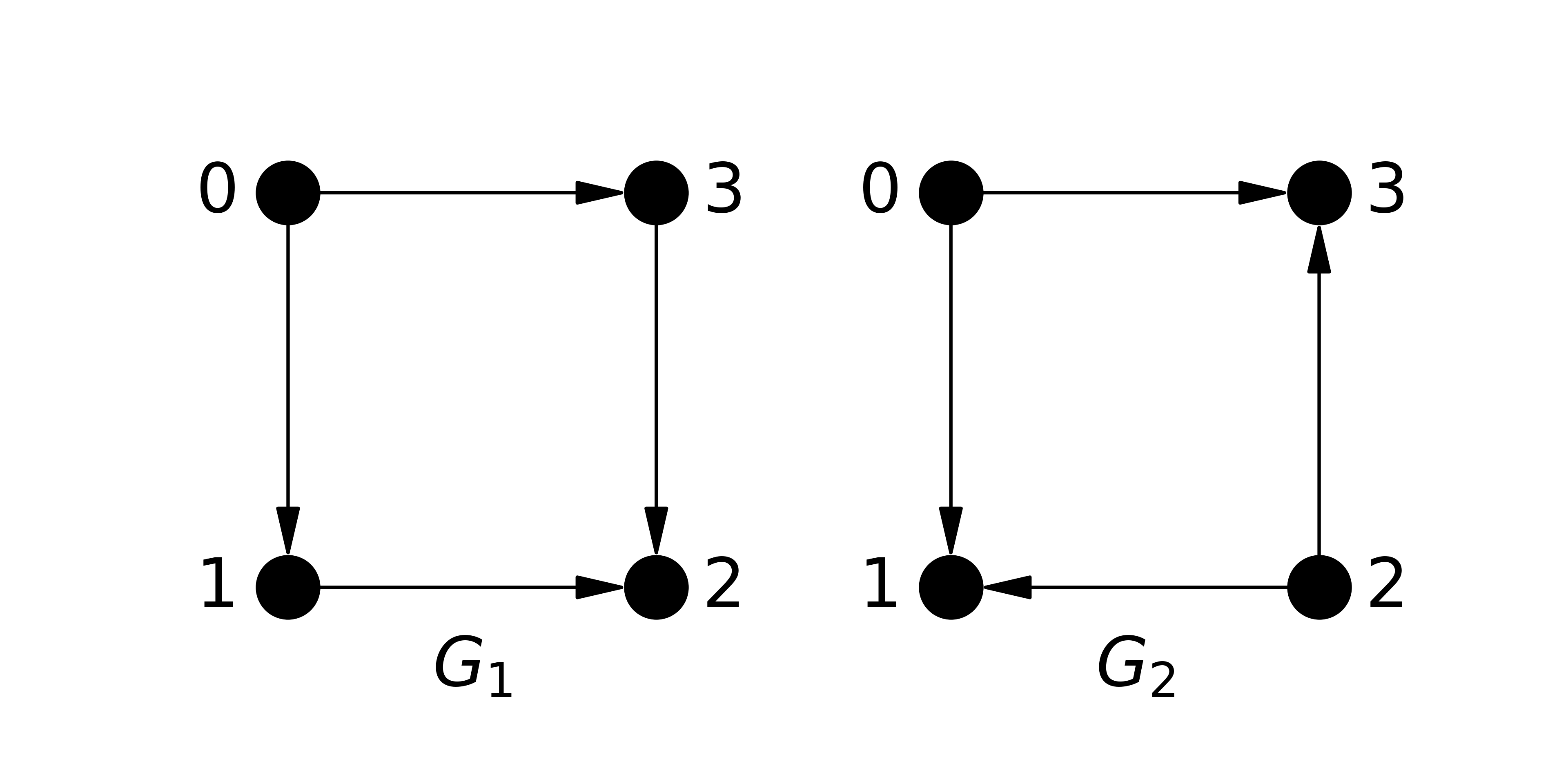}
    \caption{Directed graphs $G_1$ and $G_2$.}
    \label{fig:squares}
\end{figure}

Consider the two digraphs in Fig. \ref{fig:squares}. The path Laplacians of these digraphs are computed in Table 1 and Table 2 of \cite{wang2023persistent}. The calculation of $\Omega_2$ in $G_1$ will require the orthonormal basis as discussed above, but the calculations for $G_1$ that do not involve $\Omega_2$ not impacted as the other matrices are written with respect to orthonormal bases. Calculations for $G_2$ are also not impacted by this change. 

We recalculate the path Laplacians and spectra for $G_1$ in Table \ref{tab:path_laplacian_g1}. The only difference between the spectra here and the non-normalized version in \cite{wang2023persistent} is that here $(\lambda_1^P)_3=2$, instead of $4$, and $(\lambda_2^P)_0=2$ instead of $4$. Otherwise the spectra are the same, including the Betti numbers. We also record the calculations of the path Laplacians and spectra for $G_2$ in Table \ref{tab:path_laplacian_g2}, which are the same as in Table 2 of \cite{wang2023persistent}.

\begin{table}[]
    \centering
    \begin{tabular}{|c|c|c|c|}\hline
        $n$ & $n=0$ & $n=1$ & $n=2$\\\hline
        $\Omega_n$ & $\mathrm{span}\{e_0,e_1,e_2,e_3\} $& $\mathrm{span}\{e_{01},e_{03},e_{12},e_{32}\} $& $\mathrm{span}\{\frac{1}{\sqrt{2}}e_{032}-\frac{1}{\sqrt{2}}e_{012}\} $ \\\hline
        $B_{n+1}$ & $\begin{blockarray}{ccccc}
            & e_{01} & e_{03} & e_{12} & e_{32} \\
            \begin{block}{c(cccc)}
              e_0 & -1 & -1 & 0 & 0 \\
              e_1 & 1 & 0 & -1 & 0 \\
              e_2 & 0 & 0 & 1 & 1 \\
              e_3 & 0 & 1 & 0 & -1\\
            \end{block}
        \end{blockarray}$ & $\begin{blockarray}{cc}
            & \frac{1}{\sqrt{2}}e_{032}-\frac{1}{\sqrt{2}}e_{012} \\
            \begin{block}{c(c)}
              e_{01} & \frac{-1}{\sqrt{2}} \\
              e_{03} & \frac{1}{\sqrt{2}} \\
              e_{12} & \frac{-1}{\sqrt{2}}  \\
              e_{32} & \frac{1}{\sqrt{2}} \\
            \end{block}
        \end{blockarray}$& $1\times 0$ empty matrix\\
        $L_n^P$ & $\begin{pmatrix}
            2 & -1 & 0 & -1\\
            -1 & 2 & -1 & 0\\
            0 & -1 & 2 & -1\\
            -1 & 0 & -1 & 2
        \end{pmatrix}$
         & $\begin{pmatrix}
             \frac{5}{2} & \frac{1}{2} & \frac{-1}{2} & \frac{-1}{2}\\
             \frac{1}{2} & \frac{5}{2} & \frac{-1}{2} & \frac{-1}{2}\\
             \frac{-1}{2} & \frac{-1}{2} & \frac{5}{2} & \frac{1}{2}\\
             \frac{-1}{2} & \frac{-1}{2} & \frac{1}{2} & \frac{5}{2}
         \end{pmatrix}$ & 
         $\begin{pmatrix}
             2
         \end{pmatrix}$\\\hline
         $\mathrm{spectra}(L_n^P)$ & $\{0,2,2,4\}$ & $\{2,2,2,4\}$ & \{2\}\\\hline
         $\beta_n$ & 1& 0 & 0\\\hline
    \end{tabular}
    \caption{Path Laplacian and Spectra computation for graph $G_1$ in Fig. \ref{fig:squares}.}
    \label{tab:path_laplacian_g1}
\end{table} 

\begin{table}[]
    \centering
    \begin{tabular}{|c|c|c|c|}\hline
        $n$ & $n=0$ & $n=1$ & $n=2$\\\hline
        $\Omega_n$ & $\mathrm{span}\{e_0,e_1,e_2,e_3\} $& $\mathrm{span}\{e_{01},e_{03},e_{21},e_{23}\} $& $\{0\} $ \\\hline
        $B_{n+1}$ & $\begin{blockarray}{ccccc}
            & e_{01} & e_{03} & e_{21} & e_{23} \\
            \begin{block}{c(cccc)}
              e_0 & -1 & -1 & 0 & 0 \\
              e_1 & 1 & 0 & 1 & 0 \\
              e_2 & 0 & 0 & -1 & -1 \\
              e_3 & 0 & 1 & 0 & 1\\
            \end{block}
        \end{blockarray}$ & $4 \times 0$ empty matrix & $0 \times 0$ empty matrix\\
        $L_n^P$ & $\begin{pmatrix}
            2 & -1 & 0 & -1\\
            -1 & 2 & -1 & 0\\
            0 & -1 & 2 & -1\\
            -1 & 0 & -1 & 2
        \end{pmatrix}$
         & $\begin{pmatrix}
             2 & 1 & 1 & 0\\
             1 & 2 & 0 & 1\\
             1 & 0 & 2 & 1\\
             0 & 1 & 1 & 2
         \end{pmatrix}$ & 
         $0 \times 0$ empty matrix\\\hline
         $\mathrm{spectra}(L_n^P)$ & $\{0,2,2,4\}$ & $\{0,2,4,4\}$ & $\emptyset$\\\hline
         $\beta_n$ & 1& 1 & 0\\\hline
    \end{tabular}
    \caption{Path Laplacian and Spectra computation for graph $G_2$ in Fig. \ref{fig:squares}.}
    \label{tab:path_laplacian_g2}
\end{table} 

Notice that the path Laplacian of $G_1$ has $\beta_1=0$ while $G_2$ gives $\beta_1=1$. Thus the path homology and path Laplacian can distinguish these graphs.

Now we compute the directed flag Laplacian for these graphs. The key difference between the path complexes and the directed flag complexes on these graphs is that $2$-paths in a path complex do not correspond to $2$-simplices ($3$-cliques) in a directed flag complex. In $G_1,$ the $2$-paths $e_{012}$ and $e_{032}$ generate $\mathcal{A}_2$, and $e_{012}-e_{032}$ generates $\Omega_2$. However, in the directed flag complex of $G_2$, there are no $3$-cliques. We compute the directed flag Laplacian of $G_1$ in Table \ref{tab:flag_laplacian_g1} and of $G_2$ in Table \ref{tab:flag_laplacian_g2}. Note that $C_2=0$ for both. The resulting spectra of $G_1$ and $G_2$ are identical, although $L_1$ has slight differences between the two. The directed flag Laplacian spectra alone do not detect the difference between these two digraphs.

\begin{table}[!htbp]
    \centering
    \begin{tabular}{|c|c|c|c|}\hline
        $n$ & $n=0$ & $n=1$ \\\hline
        $C_n$ & $\mathrm{span}\{(1),(2),(3),(4)\}$ & $\mathrm{span}\{(1,2),(1,4),(2,3),(4,3)\}$ \\\hline
        $B_{n+1}$ & $\begin{blockarray}{ccccc}
            & (1,2) & (1,4) & (2,3) & (4,3) \\
            \begin{block}{c(cccc)}
              (1) & -1 & -1 & 0 & 0 \\
              (2) & 1 & 0 & -1 & 0 \\
              (3) & 0 & 0 & 1 & 1 \\
              (4) & 0 & 1 & 0 & -1\\
            \end{block}
        \end{blockarray}$ & $4\times 0$ empty matrix \\
        $L_n$ & $\begin{pmatrix}
            2 & -1 & 0 & -1\\
            -1 & 2 & -1 & 0\\
            0 & -1 & 2 & -1\\
            -1 & 0 & -1 & 2
        \end{pmatrix}$ & $\begin{pmatrix}
            2 & 1 & -1 & 0\\
            1 & 2 & 0 & -1\\
            -1 & 0 & 2 & 1\\
            0 & -1 & 1 & 2
        \end{pmatrix}$ \\\hline
        $\mathrm{spectra}(L_n)$ &$\{0,2,2,4\}$ & $\{0,2,2,4\}$ \\\hline
        $\beta_n$ & 1 & 1 \\\hline
    \end{tabular}
    \caption{Directed flag Laplacian calculations for graph $G_1$ in Fig. \ref{fig:squares}.}
    \label{tab:flag_laplacian_g1}
\end{table}

\begin{table}[!htbp]
    \centering
    \begin{tabular}{|c|c|c|c|}\hline
        $n$ & $n=0$ & $n=1$ \\\hline
        $C_n$ & $\mathrm{span}\{(1),(2),(3),(4)\}$ & $\mathrm{span}\{(1,2),(1,4),(3,2),(3,4)\}$ \\\hline
        $B_{n+1}$ & $\begin{blockarray}{ccccc}
            & (1,2) & (1,4) & (3,2) & (3,4) \\
            \begin{block}{c(cccc)}
              (1) & -1 & -1 & 0 & 0 \\
              (2) & 1 & 0 & 1 & 0 \\
              (3) & 0 & 0 & -1 & -1 \\
              (4) & 0 & 1 & 0 & 1\\
            \end{block}
        \end{blockarray}$ & $4\times 0$ empty matrix \\
        $L_n$ & $\begin{pmatrix}
            2 & -1 & 0 & -1\\
            -1 & 2 & -1 & 0\\
            0 & -1 & 2 & -1\\
            -1 & 0 & -1 & 2
        \end{pmatrix}$ & $\begin{pmatrix}
            2 & 1 & 1 & 0\\
            1 & 2 & 0 & 1\\
            1 & 0 & 2 & 1\\
            0 & 1 & 1 & 2
        \end{pmatrix}$ \\\hline
        $\mathrm{spectra}(L_n)$ &$\{0,2,2,4\}$ & $\{0,2,2,4\}$ \\\hline
        $\beta_n$ & 1 & 1 \\\hline
    \end{tabular}
    \caption{Directed flag Laplacian calculations for graph $G_2$ in Fig. \ref{fig:squares}.}
    \label{tab:flag_laplacian_g2}
\end{table}

However, if we consider the eigenvectors corresponding to these eigenvalues, we are able to distinguish the two digraphs. In Table \ref{tab:eigenvectors} we show a linearly independent set of unit eigenvectors corresponding to the eigenvalues of $L_1^j$, the $L_1$ matrix for the graphs $G_j$, where $j=1,2$. Ordering the eigenvalues as $\lambda_0=0, \lambda_1 =2, \lambda_2=2,$ and $\lambda_3=4$, we label the eigenvectors $v_i^j$ of $L_1$ for the graph $G_j$ corresponding to the eigenvalue $\lambda_i$, for $j=1,2$ and $i=0,1,2,3$. Critically, the eigenvectors $v_i^1$ are either not eigenvectors of $L_1^2$ or, if they are, they may correspond to different eigenvalues. The same is true of the eigenvectors of $L_1^2$. For example, $v_0^1$ corresponds to $\lambda_0=0$, but $L_1^2(v_0^1)= 2 v_0^1$, i.e. it is an eigenvector with eigenvalue $2$, not $0$. This implies that two Laplacians have distinct kernels, which detects the difference between the graphs. We can also compute that $v_1^1, v_2^1$ are not eigenvectors of $L_1^2$ and $v_1^2, v_2^2$ are not eigenvectors of $L_1^1$. This is another way to distinguish the two graphs.  

\begin{table}[!htbp]
    \centering
    \begin{tabular}{|c|c|c|c|c|}\hline
            i   & 0 & 1 & 2 & 3\\\hline
            $\lambda_i$ & 0 & 2 & 2 & 4\\\hline
         $v_i^1$ & $\frac{1}{2}\begin{pmatrix}1\\-1\\1\\-1\end{pmatrix}$  &  $\begin{pmatrix}
             \frac{1}{\sqrt{6}} \\ \frac{1}{\sqrt{3}} \\ \frac{1}{\sqrt{3}} \\ \frac{1}{\sqrt{6}}
         \end{pmatrix}$  & $\begin{pmatrix}
             -\frac{1}{\sqrt{3}} \\ \frac{1}{\sqrt{6}} \\ \frac{1}{\sqrt{6}} \\ -\frac{1}{\sqrt{3}}
         \end{pmatrix}$  & $\frac{1}{2}\begin{pmatrix}-1\\-1\\1\\1\end{pmatrix}$  \\\hline
         $v_i^2$ & $\frac{1}{2}\begin{pmatrix}-1\\1\\1\\-1\end{pmatrix}$  &  $\begin{pmatrix}
             -\frac{1}{\sqrt{6}} \\ -\frac{1}{\sqrt{3}} \\ \frac{1}{\sqrt{3}} \\ \frac{1}{\sqrt{6}}
         \end{pmatrix}$  & $\begin{pmatrix}
             \frac{1}{\sqrt{3}} \\ -\frac{1}{\sqrt{6}} \\ \frac{1}{\sqrt{6}} \\ -\frac{1}{\sqrt{3}}
         \end{pmatrix}$  & $\frac{1}{2}\begin{pmatrix}1\\1\\1\\1\end{pmatrix}$\\\hline
    \end{tabular}
    \caption{Linearly independent unit eigenvectors of the directed flag Laplacian $L_1$ calculated in Tables \ref{tab:flag_laplacian_g1} and \ref{tab:flag_laplacian_g2} for graphs $G_1$ and $G_2$ of Fig. \ref{fig:squares} corresponding to the eigenvalues $\lambda\in\{0,2,2,4\}$. Eigenvectors $v_i^j$ are from the graph $G_j$.}
    \label{tab:eigenvectors}
\end{table}

We have shown in Tables \ref{tab:path_laplacian_g1}, \ref{tab:path_laplacian_g2}, \ref{tab:flag_laplacian_g1}, 
 and \ref{tab:flag_laplacian_g2}  that the most natural path complex, the one composed of all paths in a directed graph $G_2$, may have different spectra than the directed flag complex on the same directed graph. We can also see that they may agree, as in $G_1$. Therefore, these are two distinct Laplacian theories for directed graphs, and one may be more appropriate for analyzing a given situation than the other. Considering only the eigenvalues, it may appear that the path Laplacian can sometimes distinguish graphs that the directed flag Laplacian cannot, but analyzing the eigenvectors of the directed flag Laplacian may provide enough insight to still detect the difference between the graphs. Their advantages in specific applications are to be further studied in the future.

\subsection{Relationship with the hyperdigraph Laplacian}

It is also interesting to analyze the relation of the proposed method with the 
persistent hyperdigraph Laplacian   \cite{chen2023persistent}. 
We now describe the general idea of the hyperdigraph Laplacian as given in \cite{chen2023persistent}, so that one can better understand the proposed directed flag Laplacian.

Set $V$ to be a nonempty finite ordered set. A directed $p$-hyperedge is a sequence $v_0v_1\cdots v_p$ of distinct elements in $V$. A hyperdigraph is a pair $\vec{\mathcal{H}}=(V,\vec{E})$, where $\vec{E}$ is a set of hyperedges. With appropriate restrictions, one can define a homology theory on a subset of the span of $p$-hyperedges, which is called the hyperdigraph homology $H_p(\vec{\mathcal{H}};\mathbb{R})$. By taking coefficients in $\mathbb{R}$, one can define an inner product, obtain a dual to the boundary operator, and craft a Laplacian $\Delta_p^{\vec{\mathcal{H}}}$ for a hyperdigraph $\vec{\mathcal{H}}$. It has desirable properties similar to other topological Laplacians, e.g., it is self-adjoint, positive semidefinite, and $\ker\Delta_p^{\vec{\mathcal{H}}}\cong H_p(\vec{\mathcal{H}};\mathbb{R})$, so this Laplacian also recovers the Betti numbers of the homology theory.

\begin{figure}[!ht]
    \centering
    \includegraphics[width=0.6\textwidth]{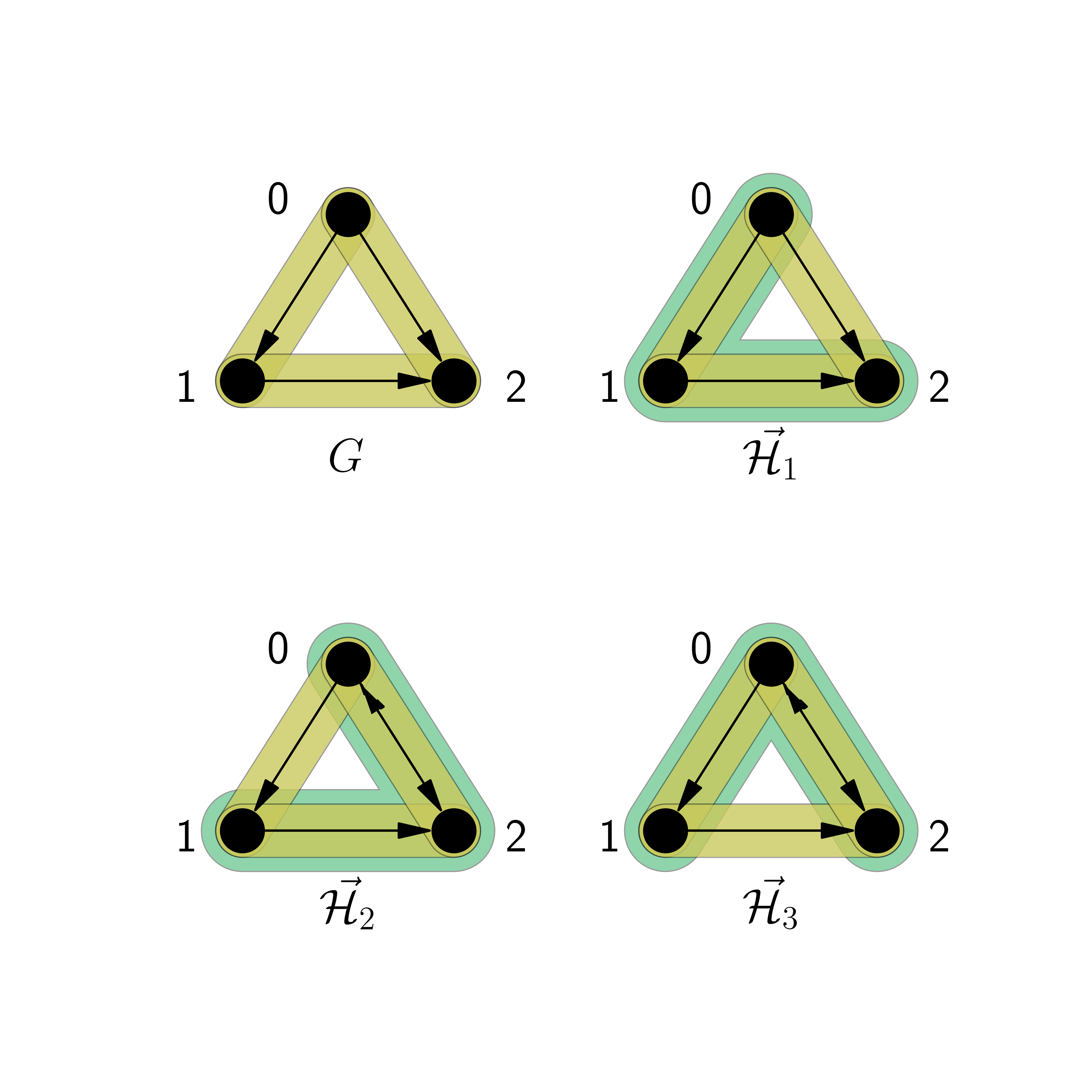}
    \caption{A digraph $G$ with three vertices and edges, and the three possible hyperdigraphs $\vec{\mathcal{H}}_1,\vec{\mathcal{H}}_2$ and $,\vec{\mathcal{H}}_3$ that have this underlying graph, but with the added hyperedges $(012),(120),$ and $(201)$, respectively. }
    \label{fig:hyperdigraph_comparison}
\end{figure}

In Fig. \ref{fig:hyperdigraph_comparison}, we show multiply hyperdigraphs with the same underlying digraph $G$. The digraph $G$ has vertices $(0),(1),$ and $(2)$, edges $(0,1),(1,2)$, and $(0,2)$. Then in $\vec{\mathcal{H}}_1$ we add the $2$-hyperedge $(012)$, in $\vec{\mathcal{H}}_2$ we add $(120)$, and in $\vec{\mathcal{H}}_3$ we add $(201)$. Hyperdigraphs allow us to add these arbitrary $2$-hyperedges, while a directed flag complex would only have the $2$-clique $(0,1,2)$. Note that $\vec{\mathcal{H}}_1$ is a hyperdigraph corresponding to the directed flag complex on $G$. This puts the extreme flexibility of the hyperdigraph on display, because adding any of the possible $2$-hyperedges on the three vertices produces a valid hyperdigraph.

Hyperdigraphs generalize abstract simplicial complexes on a finite set, so any directed flag complex can be formulated as an equivalent hyperdigraph. A directed flag complex is entirely determined by the directed edges and vertices. Unlike with path complexes, the extreme flexibility of hyperdigraphs means there is no natural hyperdigraph on a vertex set that is solely determined by vertices and directed edges. This flexibility is one of the core benefits of hyperdigraphs, since its relations could depend on higher dimensional information than pairwise connections. One could assign a $p$-hyperedge for each path of length $p$ with distinct vertices, but this would become a path complex on the directed graph.

The additional structure enforced on the directed flag complex means that it tells us about how highly connected the graph is. Because a simplex is only added when a directed clique is present, this may reduce the number of higher dimensional simplices one uses in studying the dataset. We suspect this is one of the reasons the present work is able to continue the filtration in Fig. \ref{fig:1a99_spectra} to a cutoff of $8$\r{A}, whereas the approach in Ref.  \cite{chen2023persistent} could only reach up to $4$\r{A}. The other possible reason is due to implementation factors, such as our use of  FLAGSER  to build the underlying complex, which is very fast with a mostly C++ implementation.
A more comprehensive comparison of these two  formalisms is needed but is out of the scope of the present work. In particular, there is a pressing need to develop more efficient computational algorithms for PTLs.

\section{Application}\label{sec:application}

\begin{figure}
    \centering
    \includegraphics[width=0.6\textwidth]{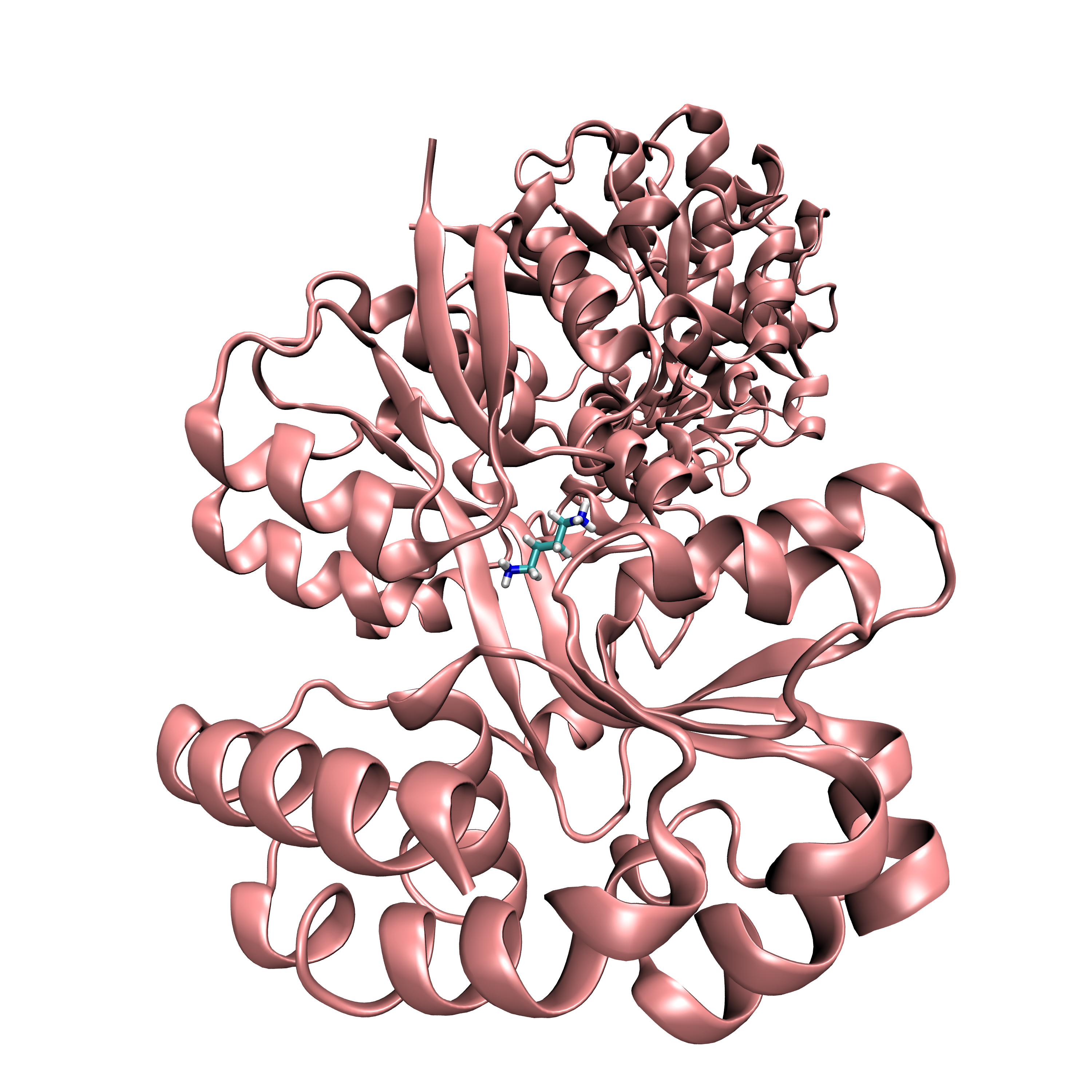}
    \caption{Protein-Ligand complex with PDBID 1a99.}
    \label{fig:1a99}
\end{figure}

\begin{figure}
    \centering
    \includegraphics[width=0.8\textwidth]{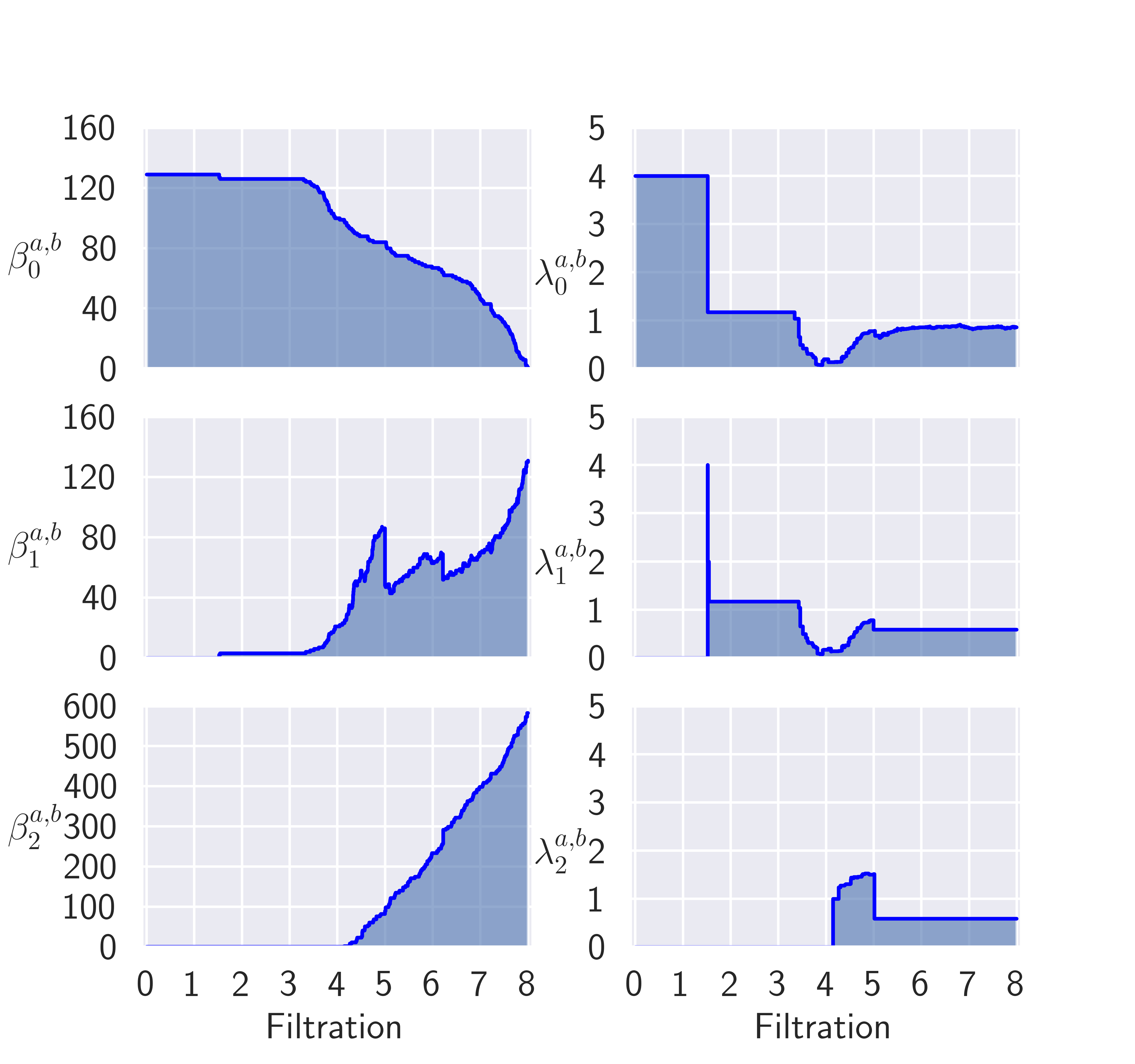}
    \caption{Spectra of protein-ligand complex with PDB ID 1a99, with restrictions on included atoms and filtered edges. Note that $\beta_2^{a,b}$ uses a larger vertical axis scale than $\beta_0^{a,b}$ and $\beta_1^{a,b}$. }
    \label{fig:1a99_spectra}
\end{figure}

In this section, we demonstrate the PDFL by analyzing   protein-ligand complex 1a99 from the Protein Data Bank (PDB) as shown in Fig. \ref{fig:1a99} for visualizing the complex. The same complex was studied in the literature \cite{chen2023persistent}. To analyze the protein-ligand binding,  we consider only C atoms of the protein within $8$ \r{A} of the ligand, and only C, N, O, and S atoms of the ligand. We use the distance between two atoms in this restricted complex as the filtration value at which an edge is added between the two atoms. We round these distances to the nearest $10^{-3}$. We add edges for the bonds within the ligand, but not for bonds within the protein. The direction of the edge is pointing toward the atom with the higher electronegativity of the pair; in the case of a C-C pair, we add both possible edges. We then compute $\lambda_k^{a,b}$ and $\beta_k^{a,b}$, where $k=0,1,2$, and $a$ is the filtration parameter, and $b$ is the next filtration parameter greater than $a$. 
We used a cutoff value of $8$\r{A} angstrom because it was sufficiently large to allow the non-harmonic spectra to stabilize, and for a clear pattern to emerge in the harmonic spectra, but the implementation would allow for reasonable computation time to determine the spectra associated with larger cutoffs in different systems.

The computed spectra are in Fig. \ref{fig:1a99_spectra}. We can make several observations about this system. First, observe that $\beta_0^{a,b}$ tends to $1$, since the final graph is fully connected. Next, observe that $\beta_2^{a,b}$ increases essentially monotonically; this is because the particular graph we produced does not have $3$-simplices. It is challenging to have any $3$-simplex because we do not add edges for within-protein atoms, and so there can only be $1$ protein atom in any clique. To make a $3$-simplex, the lone protein atom would need to have edges with $3$ ligand atoms, that are all ``bonded'' to each other. This occurs rarely in nature. 

Another notable feature of the spectra is that all are essentially constant until the filtration value is $3.$ From there on, $\beta_0^{a,b}$ decreases steadily, which reflects that we added edges at every filtration value, and it is not common for many edges to have the same filtration value, which comes from the distance between atoms, rounded to the nearest $10^{-3}$. We can see that $\beta_1^{a,b}$ rises rapidly from a filtration of $3$ to $5$, and then sharply decreases at $5$ before continuing to increase, with some other perturbations. We can see that $\lambda_0^{a,b}$, $\lambda_1^{a,b}$, and $\lambda_2^{a,b}$ all change suddenly for filtration values around $3$ and $5$, with a sharp change at $5$, but what is interesting is that they stabilize thereafter even as $\beta_1^{a,b}$ continues to change. This suggests that there are few significant interactions between the ligand molecule and protein atoms further than $5$\r{A}.

\section{Conclusion}\label{sec:conclusion}

Recently, persistent topological Laplacians (PTLs) have been proposed to overcome many drawbacks of classical persistent homology in topological data analysis (TDA), giving rise to a powerful new way to study the filtrations of simplicial complexes and other structures. One type of abstract simplicial complexes, directed flag complexes, commonly arises in studies of networks, including small world networks,  neuron networks, and gene regulation networks. Therefore,  a theory of persistent directed flag Laplacians (PDFL) may provide additional insight into challenging problems. This work builds such a theory and shows that it can provide a more nuanced representation of persistent changes in directed flag complexes than persistent Betti numbers alone can. Comparison is given to earlier persistent path Laplacians and persistent hyperdigraph Laplacians. We also demonstrate how PDFL can help analyze networks arising from real-world data, including a protein-ligand system. These examples illustrate the application potential of proposed PDFL to practice problems.

\section*{Code availability}

The codes are available at  \href{https://github.com/bdjones13/flagser-laplacian/}{https://github.com/bdjones13/flagser-laplacian/}.

\section*{Acknowledgments}
This work was supported in part by NIH grants R01GM126189, R01AI164266, and R35GM148196, National Science Foundation grants DMS2052983, DMS-1761320, and IIS-1900473, NASA  grant 80NSSC21M0023,   Michigan State University Research Foundation, and  Bristol-Myers Squibb  65109.  The authors thank Dong Chen for technical assistance and Xiaoqi Wei for proofreading.  

\bibliographystyle{abbrv}

\bibliography{main_v7}

\end{document}